\documentclass[11pt,a4paper]{article}
\usepackage{amsmath,amsfonts,amssymb,amsthm,latexsym}
\usepackage{graphics,epsfig,color,enumerate}
\usepackage{setspace,mathtools,enumitem,dsfont,verbatim}
\usepackage{caption}
\usepackage{subcaption}
\usepackage[english]{babel}
\makeatletter
\def\@captype{figure}
\makeatother

\newtheorem{theorem}{Theorem}[section]

\newtheorem{remark}[theorem]{Remark}

\numberwithin{equation}{section}

\newcommand{\R}{\mathbb{R}}

\newcommand{\Z}{\mathbb{Z}}
\newcommand{\N}{\mathbb{N}}

\newcommand{\cG}{\mathcal{G}}
\newcommand{\cP}{\mathcal{P}}

\newcommand{\bG}{\mathbf{G}}
\newcommand{\bGamma}{\mathbf{\Gamma}}
\newcommand{\vC}{\vec{C}}
\newcommand{\vt}{\vec{\theta}}
\newcommand{\Pmic}{P_{\mathrm{mic}}}
\newcommand{\Pcan}{P_{\mathrm{can}}}

\newcommand{\be}{\begin{equation}}
\newcommand{\ee}{\end{equation}}

\title{Ensemble nonequivalence in random graphs\\ 
with modular structure}

\author{

Diego Garlaschelli
\footnotemark[1]
\\

Frank den Hollander
\footnotemark[2]
\\

Andrea Roccaverde
\footnotemark[12]
}

\footnotetext[1]{
Lorentz Institute for Theoretical Physics, Leiden University, P.O.\  Box 9504, 
2300 RA Leiden, The Netherlands
}

\footnotetext[2]{
Mathematical Institute, Leiden University, P.O.\ Box 9512,
2300 RA Leiden, The Netherlands
}

\date{\today}

\begin{document}

\maketitle 

\begin{abstract}
Breaking of equivalence between the microcanonical ensemble and the canonical ensemble,
describing a large system subject to hard and soft constraints, respectively, was recently shown 
to occur in large random graphs. Hard constraints must be met by every graph, soft constraints 
must be met only on average, subject to maximal entropy. In Squartini, de Mol, den Hollander 
and Garlaschelli (2015) it was shown that ensembles of random graphs are nonequivalent when
the degrees of the nodes are constrained, in the sense of a non-zero limiting specific relative 
entropy as the number of nodes diverges. In that paper, the nodes were placed either on a 
single layer (uni-partite graphs) or on two layers (bi-partite graphs). In the present paper we 
consider an arbitrary number of intra-connected and inter-connected layers, thus allowing for 
modular graphs with a multi-partite, multiplex, time-varying, block-model or community structure. We give 
a full classification of ensemble equivalence in the sparse regime, proving that breakdown 
occurs as soon as the number of local constraints (i.e., the number of constrained degrees) 
is extensive in the number of nodes, irrespective of the layer structure. In addition, we derive 
an explicit formula for the specific relative entropy and provide an interpretation of this formula in terms 
of Poissonisation of the degrees.



\medskip\noindent
{\it Acknowledgment.} 
DG and AR are supported by EU-project 317532-MULTIPLEX. FdH and AR are supported by NWO 
Gravitation Grant 024.002.003--NETWORKS.
\end{abstract}

\newpage


\section{Introduction and main results}
\label{S1}


\subsection{Background and outline}
\label{S1.1}

For systems with many interacting components a detailed microscopic description is 
infeasible and must be replaced by a probabilistic description, where the system is 
assumed to be a random sample drawn from a set of allowed microscopic configurations 
that are consistent with a set of known macroscopic properties, referred to as \emph{constraints}. 
Statistical physics deals with the definition of the appropriate probability distribution 
over the set of microscopic configurations and with the calculation of the resulting
macroscopic properties of the system. The three main choices of probability distribution 
are: (1) the \emph{microcanonical ensemble}, where the constraints are \emph{hard} 
(i.e., are satisfied by each individual configuration); (2) the \emph{canonical ensemble}, 
where the constraints are \emph{soft} (i.e., hold as ensemble averages, while individual 
configurations may violate the constraints); (3) the \emph{grandcanonical ensemble}, 
where also the number of components is considered as a soft constraint.

For systems that are large but finite, the three ensembles are obviously different and, in fact, 
represent different physical situations: (1) the microcanonical ensemble models completely 
isolated systems (where both the energy and the number of particles are ``hard''); (2) the 
canonical ensemble models closed systems in thermal equilibrium with a heat bath (where 
the energy is  ``soft'' and the number of particles is ``hard''); (3) the grandcanonical ensemble 
models open systems in thermal and chemical equilibrium (where both the energy and the 
number of particles are ``soft''). However, in the limit as the number of particles diverges, the 
three ensembles are traditionally \emph{assumed} to become equivalent as a result of the 
expected vanishing of the fluctuations of the soft constraints, i.e., the soft constraints are expected 
to become asymptotically hard. This assumption of \emph{ensemble equivalence}, which dates
back to Gibbs~\cite{Gibbs1902}, has been verified in traditional models of physical systems with 
short-range interactions and a finite number of constraints, but it does \emph{not} hold in general. 
Nonetheless, equivalence is considered to be one of the pillars of statistical physics and underlies 
many of the results that contribute to our current understanding of large real-world systems.

Despite the fact that many textbooks still convey the message that ensemble equivalence holds 
for all systems, as some sort of universal asymptotic property, over the last decades various 
examples have been found for which it breaks down. These examples range from astrophysical 
processes \cite{LyndenBell1968}, \cite{Thirring1970}, \cite{Thirring1971}, \cite{LyndenBell1999},
\cite{Chavanis2003}, quantum phase separation \cite{Blume1971}, \cite{Barre2001}, \cite{Ellis2004}, 
nuclear fragmentation \cite{DAgostino2000}, and fluid turbulence \cite{Ellis2000}, \cite{Ellis2002}. 
Across these examples, the signatures of ensemble nonequivalence differ, which calls for a 
rigorous mathematical definition of ensemble (non)equivalence: (i) \emph{thermodynamic equivalence} 
refers to the existence of an invertible Legendre transform between the microcanonical entropy 
and canonical free energy \cite{Ellis2004}; (ii) \emph{macrostate equivalence} refers to the equivalence 
of the canonical and microcanonical sets of equilibrium values of macroscopic properties 
\cite{Touchette2004}; (iii) \emph{measure equivalence} refers to the asymptotic equivalence 
of the microcanonical and canonical probability distributions in the thermodynamic limit, i.e., the 
vanishing of their specific relative entropy \cite{touchette2014}. The latter reference reviews the 
three definitions and shows that, under certain hypotheses, they are identical. 

In the present paper we focus on the equivalence between microcanonical and canonical ensembles, 
although nonequivalence can in general involve the grandcanonical ensemble as well~\cite{ziff}. While 
there is consensus that nonequivalence occurs when the microcanonical specific entropy is non-concave 
as a function of the energy density in the thermodynamic limit,  the classification of the 
physical mechanisms at the \emph{origin} of nonequivalence is still open. In most 
of the models studied in the literature, nonequivalence appears to be associated with the 
non-additivity of the energy of the subparts of the system or with phase transitions \cite{Ruffo}, 
\cite{Ruffo2}, \cite{touchette2014}. A possible and natural mechanism for non-additivity is the 
presence of \emph{long-range interactions}. Similarly, phase transitions are naturally associated 
with long-range order. These ``standard mechanisms'' for ensemble nonequivalence have 
been documented also in the study of random graphs. In \cite{Barre2007}, a Potts model on a 
random regular graph is studied in both the microcanonical and canonical ensemble, where 
the microscopic configurations are the spin configurations (not the configurations of the network 
itself). It is found that the long-range nature of random connections, which makes the model 
non-additive and the microcanonical entropy non-concave, ultimately results in ensemble 
nonequivalence. In~\cite{RS}, ~\cite{RS3}, ~\cite{RS2} and~\cite{chatterjee}, random networks 
with given densities of edges and triangles are considered, and phase transitions characterised 
by jumps in these densities are found, with an associated breaking of ensemble equivalence 
(where the microscopic configurations are network configurations).

Recently, the study of certain classes of uni-partite and bi-partite random graphs \cite{SdMdHG15},
\cite{GdHR15} has shown that ensemble nonequivalence can manifest itself via an additional, 
novel mechanism, unrelated to non-additivity or phase transitions: namely, the presence of an 
\emph{extensive} number of local topological constraints, i.e., the degrees and/or the strengths 
(for weighted graphs) of all nodes.\footnote{While in binary (i.e., simple) graphs the \emph{degree} 
of a node is defined as the number of edges incident to that node, in weighted graphs (i.e., 
graphs where edges can carry weights) the \emph{strength} of a node is defined as the total 
weight of all edges incident to that node. In this paper, we focus on binary graphs only.} This 
finding explains previously documented signatures of nonequivalence in random graphs with local 
constraints, such as a finite difference between the microcanonical and canonical entropy densities 
\cite{ginestra} and the non-vanishing of the relative fluctuations of the constraints \cite{myunbiased}.
How generally this result holds beyond the specific uni-partite and bi-partite cases considered so far 
remains an open question, on which we focus in the present paper. By considering a much more 
general class of random graphs with a variable number of constraints, we confirm that the presence 
of an extensive number of local topological constraints breaks ensemble equivalence, even in the 
absence of phase transitions or non-additivity. 

The remainder of our paper is organised as follows. In Section~\ref{S1.2} we give the definition 
of measure equivalence and, following \cite{SdMdHG15}, show that it translates into a simple 
pointwise criterion for the large deviation properties of the microcanonical and canonical 
probabilities. In Section~\ref{S1.3} we introduce our main theorems in pedagogical order, starting 
from the characterisation of nonequivalence in the simple cases of uni-partite and bi-partite 
graphs already explored in~\cite{SdMdHG15}, and subsequently moving on to a very general 
class of graphs with arbitrary multilayer structure and tunable intra-layer and inter-layer connectivity.
Our main theorems, which (mostly) concern the \emph{sparse regime}, not only characterise 
nonequivalence \emph{qualitatively}, they also provide a \emph{quantitative} formula for the 
specific relative entropy. In Section~\ref{S2} we discuss various important implications of our 
results, describing properties that are fully general but also focussing on several special cases 
of empirical relevance. In addition, we provide an interpretation of the specific relative entropy 
formula in terms of Poissonisation of the degrees. We also discuss the implications of our results 
for the study of several empirically relevant classes of ``modular'' networks that have recently 
attracted interest in the literature, such as networks with a so-called multi-partite, multiplex 
\cite{multiplex}, time-varying \cite{timevarying}, block-model \cite{simpleBM}, \cite{degreecorrectedBM} 
or community structure \cite{community_santo}, \cite{community_porter}. In Section~\ref{S3}, 
finally, we provide the proofs of our theorems.

In future work we will address the \emph{dense regime}, which requires the use of \emph{graphons}.
In that regime we expect nonequivalence to persist, and in some cases become even more pronounced.


\subsection{Microcanonical ensemble, canonical ensemble, relative entropy}
\label{S1.2}

For $n \in \N$, let $\cG_n$ denote the set of all simple undirected graphs with $n$ nodes. Let 
$\cG^\sharp_n\subseteq\cG_n$ be some non-empty subset of $\cG_n$, to be specified later. 
Informally, the restriction from $\cG_n$ to $\cG^\sharp_n$ allows us to forbid the presence of 
certain links, in such a way that the $n$ nodes are effectively partitioned into $M \in \N$ groups 
of nodes (or ``layers'') of sizes $n_1,\ldots,n_M$ with $\sum_{i=1}^{M} n_i =n$. This restriction 
can be made explicit and rigorous through the definition of a superstructure, which we call the 
\emph{master graph}, that will be introduced later. A given choice of $\cG^\sharp_n$ corresponds 
to the selection of a specific class of \emph{multilayer} graphs with desired intra-layer and 
inter-layer connectivity, such as graphs with a multipartite, multiplex, time-varying, block-model 
or community structure. In the simplest case, $\cG^\sharp_n=\cG_n$, which reduces to the 
ordinary choice of  uni-partite (single-layer) graphs. This example, along with various more 
complicated examples, is considered explicitly later on.

In general, any graph $\bG\in\cG^\sharp_n$ can be represented as an $n \times n$ matrix 
with elements 
\be
g_{i,j}(\bG) =
\begin{cases}
1\qquad \mbox{if there is a link between node\ } i \mbox{\ and node\ } j,\\ 
0 \qquad \mbox{otherwise.}
\end{cases}
\ee
Let $\vC$ denote a vector-valued function on $\cG^\sharp_n$. Given a specific value $\vC^*$, 
which we assume to be \emph{graphic}, i.e., realisable by at least one graph in $\cG^\sharp_n$, 
the \emph{microcanonical probability distribution} on $\cG^\sharp_n$ with \emph{hard constraint} 
$\vC^*$ is defined as
\begin{equation}
\Pmic(\bG) =
\left\{
\begin{array}{ll} 
1/\Omega_{\vC^*}, \quad & \text{if } \vC(\bG) = \vC^*, \\ 
0, & \text{else},
\end{array}
\right.
\label{eq:PM}
\end{equation}
where 
\begin{equation}
\Omega_{\vC^*} = | \{\bG \in \cG^\sharp_n\colon\, \vC(\bG) = \vC^* \} |>0
\end{equation}
is the number of graphs that realise $\vC^*$. The \emph{canonical probability distribution} 
$\Pcan(\bG)$ on $\cG^\sharp_n$ is defined as the solution of the maximisation of the 
\emph{entropy} 
\begin{equation}
S_n(\Pcan) = - \sum_{\bG \in \cG^\sharp_n} \Pcan(\bG) \ln \Pcan(\bG)
\end{equation}
subject to the \emph{soft constraint} $\langle \vC \rangle  = \vC^*$, where $\langle \cdot 
\rangle$ denotes the average w.r.t.\ $\Pcan$, and subject to the normalisation condition 
$\sum_{\bG \in \cG^\sharp_n} \Pcan(\bG) = 1$. This gives
\begin{equation}
\Pcan(\bG) = \frac{\exp[-H(\bG,\vt^*)]}{Z(\vt^*)},
\label{eq:PC}
\end{equation}
where 
\begin{equation}
H(\bG, \vt) = \vt \cdot \vC(\bG)
\label{eq:H}
\end{equation}
is the \emph{Hamiltonian} (or \emph{energy}) and
\be
Z(\vt\,) = \sum_{\bG \in \cG^\sharp_n} \exp[-H(\bG, \vt\,)]
\ee
is the \emph{partition function}. Note that in \eqref{eq:PC} the parameter $\vt$ must be set 
to the particular value $\vt^*$ that realises $\langle \vC \rangle  = \vC^*$. This value also 
maximises the likelihood of the model, given the data \cite{mylikelihood}.

It is worth mentioning that, in the social network analysis literature \cite{WF}, maximum-entropy 
canonical ensembles of graphs are traditionally known under the name of Exponential Random 
Graphs (ERGs). Indeed, many of the examples of canonical graph ensembles that we will consider 
in this paper, or variants thereof, have been studied previously as ERG models of social networks. 
Recently, ERGs have also entered the physics literature \cite{ginestra}, \cite{gin_cavity}, \cite{gin_hierarchically}, 
\cite{ParkNewman}, \cite{mymethod}, \cite{myunbiased}, \cite{myenhanced}, \cite{myX} ,\cite{myY}, 
\cite{degreecorrectedBM}, \cite{fronczakBM}, \cite{entropyBM}, 
\cite{bianconiMultiplex} because of the wide applicability of techniques from statistical 
physics for the calculation of canonical partition functions. We will refer more extensively to 
these models, and to the empirical situations for which they have been proposed, in 
Section~\ref{sec:empirical}. Apart for a few exceptions \cite{ginestra}, \cite{entropyBM}, 
\cite{SdMdHG15}, these previous studies have not addressed the problem of ensemble 
(non)equivalence of ERGs. The aim of the present paper is to do so exhaustively, and in 
a mathematically rigorous way, via the following definitions.

The \emph{relative entropy} of $\Pmic$ w.r.t.\ $\Pcan$ is
\begin{equation}
S_n(\Pmic \mid \Pcan) 
= \sum_{\bG \in \cG^\sharp_n} \Pmic(\bG) \ln \frac{\Pmic(\bG)}{\Pcan(\bG)},
\label{eq:KL1}
\end{equation}
and the \emph{specific relative entropy} is \begin{equation}
s_n = n^{-1}\,S_n(\Pmic \mid \Pcan).
\label{eq:sn}
\end{equation}
Following \cite{touchette2014}, \cite{SdMdHG15}, we say that the two ensembles are 
measure equivalent if and only if their specific relative entropy vanishes in the 
\emph{thermodynamic limit} $n\to\infty$, i.e., 
\begin{equation}
s_{\infty} = \lim_{n \to \infty} n^{-1}\,S_n(\Pmic \mid \Pcan) = 0.
\label{eq:criterion1}
\end{equation}
It should be noted that, for a given choice of $\cG^\sharp_n$ and $\vC$, there may be 
different ways to realise the thermodynamic limit, corresponding to different ways in 
which  the numbers $\{n_i\}_{i=1}^M$ of nodes inside the $M$ layers grow relatively 
to each other. So, \eqref{eq:criterion1} implicitly requires an underlying \emph{specific 
definition of the thermodynamic limit}. Explicit examples will be considered in each case 
separately, and certain different realisations of the thermodynamic limit will indeed be 
seen to lead to different results. With this in mind, we suppress the $n$-dependence 
from our notation of quantities like $\bG$, $\vC$, $\vC^*$, $\Pmic$, $\Pcan$, $H$, $Z$. 
When letting $n\to\infty$ it will be understood that $\bG \in \cG^\sharp_n$ always.

Before considering specific cases, we recall an important observation made in \cite{SdMdHG15}. 
The definition of $H(\bG,\vt\,)$ ensures that, for any $\bG_1,\bG_2\in\cG^\sharp_n$, 
$\Pcan(\bG_1)=\Pcan(\bG_2)$ whenever $\vC(\bG_1)=\vC(\bG_2)$ (i.e., the canonical 
probability is the same for all graphs having the same value of the constraint). We may 
therefore rewrite \eqref{eq:KL1} as
\begin{equation}
S_n(\Pmic \mid \Pcan) = \ln \frac{\Pmic(\bG^*)}{\Pcan(\bG^*)},
\label{eq:KL2}
\end{equation}
where $\bG^*$ is \emph{any} graph in $\cG^\sharp_n$ such that $\vC(\bG^*) =\vC^*$ 
(recall that we have assumed that $\vC^*$ is realisable by at least one graph in $\cG^\sharp_n$). 
The condition for equivalence in \eqref{eq:criterion1} then becomes  
\begin{equation}
\lim_{n \to \infty} n^{-1}\, \big[\ln {\Pmic(\bG^*)} - \ln{\Pcan(\bG^*)} \big] = 0,
\label{eq:criterion3}
\end{equation}
which shows that the breaking of ensemble equivalence coincides with $\Pmic(\bG^*)$ and 
$\Pcan(\bG^*)$ having different large deviation behaviour. Importantly, this condition is entirely 
local, i.e., it involves the microcanonical and canonical probabilities of a \emph{single} configuration 
$\bG^*$ realising the hard constraint. Apart from its theoretical importance, this fact greatly 
simplifies mathematical calculations. Note that \eqref{eq:criterion3}, like \eqref{eq:criterion1}, 
implicitly requires a specific definition of the thermodynamic limit. For a given choice of 
$\cG^\sharp_n$ and $\vC$, different definitions of the thermodynamic limit may result either 
in ensemble equivalence or in ensemble nonequivalence.


\subsection{Main Theorems}
\label{S1.3}

Most of the constraints that will be considered below are \emph{extensive} in the 
number of nodes.


\subsubsection{Single layer: uni-partite graphs}
\label{S1.3.1}

The first class of random graphs we consider is specified by $M=1$ and $\cG^\sharp_n
=\cG_n$. This choice corresponds to the class of (simple and undirected) \emph{uni-partite 
graphs}, where links are allowed between each pair of nodes. We can think of these graphs 
as consisting of a single layer of nodes, inside which all links are allowed. Note that in this 
simple case the thermodynamic limit $n\to\infty$ can be realised in a unique way, which makes 
\eqref{eq:criterion1} and \eqref{eq:criterion3} already well-defined. 

\paragraph{Constraints on the degree sequence.}
For a uni-partite graph $\bG\in\cG_n$, the degree sequence is defined as $\vec{k}(\bG) 
= (k_i(\bG))_{i=1}^n$ with $k_i(\bG)=\sum_{j\ne i}g_{i,j}(\bG)$. In what follows we 
constrain the degree sequence to a \emph{specific value} $\vec{k}^*$, which (in accordance 
with our aforementioned general prescription for $\vC^*$) we assume to be \emph{graphical}, 
i.e., there is at least one graph with degree sequence $\vec{k}^*$. The constraints are therefore
\be
\label{degcon}
\vC^* = \vec{k}^*= (k_i^*)_{i=1}^n \in \N_0^n,
\ee 
where $\N_0 = \N \cup \{0\}$ with $\N=\{1,2,\ldots\}$. This class is also known as the 
\emph{configuration model} (\cite{bender1978}, \cite{bollobas1980}, \cite{molloyreed}, 
\cite{newmanstrogatz}, \cite{chunglu}, \cite{myunbiased}; see also \cite[Chapter 7]{vdH16}). 
In \cite{SdMdHG15} the breaking of ensemble equivalence was studied in the \emph{sparse 
regime} defined by the condition
\be
m^* = \max_{1\le i\le n} k^*_i = o(\sqrt{n}).
\label{eq:sparse1}
\ee

Let $\cP(\N_0)$ denote the set of probability distributions on $\N_0$. Let
\be
\label{empdef}
f_n = n^{-1} \sum_{i=1}^n \delta_{k^*_i} \in \cP(\N_0),
\ee
be the \emph{empirical degree distribution}, where $ \delta_{k}$ denotes the point measure at $k$. 
Suppose that there exists a degree distribution $f \in \cP(\N_0)$ such that
\be
\lim_{n\to\infty} \|f_n - f\|_{\ell^1(g)} = 0,
\label{eq:conv1}
\ee
where $g\colon\N_0 \to [0,\infty)$ is given by 
\be
\label{gdef}
g(k) = \log\left(\frac{k!}{k^k e^{-k}}\right), \qquad k \in \N_0,
\ee
and $\ell^1(g)$ is the vector space of functions $h\colon\,\Z\to\R$ with $\|h\|_{\ell^1(g)} 
= \sum_{k\in\N_0} |h(k)| g(k)<\infty$. For later use we note that
\be
g(0)=0, \qquad k \mapsto g(k) \text{ is strictly increasing}, 
\qquad g(k) = \tfrac12 \log (2\pi k) + O(k^{-1}), \quad k \to\infty.
\label{g(k)}
\ee

\begin{theorem}
\label{resultunip}
Subject to \eqref{degcon}--\eqref{eq:sparse1} and \eqref{eq:conv1}, the specific relative entropy 
equals
\be
s_{\infty} = \|f\|_{\ell^1(g)}>0.
\label{sec:unis}
\ee
\end{theorem}

\noindent
Thus, when we constrain the degrees we break the ensemble equivalence. 

\begin{remark}
\label{remAL}
{\rm It is known that $\vec{k}^*$ is graphical if and only if $\sum_{i=1}^n k_i^*$ is even 
and 
\be
\label{grcond}
\sum_{i=1}^j k_i^* \leq j(j-1) + \sum_{i=j+1}^n \min (j,k^*_i), \qquad j =1,\ldots,n-1.
\ee
In \cite{arratialiggett2005}, the case where $k_i^*$, $i\in\N$, are i.i.d.\ with probability 
distribution $f$ is considered, and it is shown that
\be
\lim_{n\to\infty}  f^{\otimes n}\Big((k_1^*,\ldots,k_n^*) \text{ is graphical} 
~\Big|~ \sum_{i=1}^n k_i^* \text{ is even}\Big) = 1
\ee
as soon as $f$ satisfies $0<\sum_{k\,\text{even}} f(k) < 1$ and $\lim_{n\to\infty} n \sum_{k \geq n} 
f(k) = 0$. (The latter condition is slightly weaker than the condition $\sum_{k\in\N_0} kf(k)
<\infty$.) In what follows we do \emph{not} require the degrees to be drawn in this manner, but 
when we let $n\to\infty$ we always implicitly assume that the limit is taken \emph{within the class of 
graphical degree sequences}. } 
\end{remark}

\begin{remark}
{\rm A different yet similar definition of sparse regime, replacing \eqref{eq:sparse1}, is given in 
van der Hofstad~\cite[Chapter 7]{vdH16}. This condition is formulated in terms of bounded 
second moment of the empirical degree distribution $f_n$ in the limit as $n\to\infty$. 
Theorem~\ref{resultunip} carries over.}
\end{remark}

\paragraph{Constraints on the total number of links only.} 
We now relax the constraints, and fix only the total number of links $L(\bG)=\frac{1}{2}
\sum_{i=1}^n k_i(\bG)$. The constraint therefore becomes
\be
\label{linkcon}
\vC^* = L^*.
\ee 
It should be note that in this case, the canonical ensemble coincides with the Erd\H{o}s-R\'enyi 
random graph model, where each pair of nodes is independently connected with the same 
probability. As shown in \cite{ginestra}, \cite{SdMdHG15}, in this case the usual result that the 
ensembles are asymptotically equivalent holds.

\begin{theorem}
\label{biplinks}
Subject to \eqref{linkcon}, the specific relative entropy equals $s_{\infty}=0$.
\end{theorem}


\subsubsection{Two layers: bi-partite graphs}
\label{S1.3.2}

The second class of random graphs we consider are \emph{bi-partite graphs}. Here $M=2$ 
and nodes are placed on two (non-overlapping) layers (say, top and bottom), and only links 
\emph{across} layers are allowed. Let $\Lambda_1$ and $\Lambda_2$ denote the sets of 
nodes in the top and bottom layer, respectively. The set of all bi-partite graphs consisting of 
$n_1=|\Lambda_1|$ nodes in the top layer and ${n_2}=|\Lambda_2|$ nodes in the bottom 
layer is denoted by $\cG^\sharp_n=\cG_{n_1,{n_2}}\subset \cG_n$. Bi-partiteness means 
that, for all $\bG\in\cG_{n_1,n_2}$, we have $g_{i,j}(\bG)=0$ if $i,j\in\Lambda_1$ or $i,j\in
\Lambda_2$.

In a bipartite graph $\bG\in\cG_{n_1,n_2}$, we define the degree sequence of the top layer 
as $\vec{k}_{1\to2}(\bG)=(k_i(\bG))_{i\in\Lambda_1}$, where $k_i(\bG)=\sum_{j\in
\Lambda_2}g_{i,j}(\bG)$. Similarly, we define the degree sequence of the bottom layer as 
$\vec{k}_{2\to1}(\bG)=(k'_i(\bG))_{i\in\Lambda_2}$, where $k'_i(\bG)=\sum_{j\in
\Lambda_1}g_{i,j}(\bG)$. The symbol $s\to t$ highlights the fact that the degree sequence 
of layer $s$ is built from links pointing from $\Lambda_s$ to $\Lambda_t$ ($s,t=1,2$). The 
degree sequences $\vec{k}_{1\to2}(\bG)$ and  $\vec{k}_{2\to1}(\bG)$ are related by the 
condition that they both add up to the total number of links $L(\bG)$:
\be 
L(\bG)= \sum_{i\in\Lambda_1} k_i(\bG) = \sum_{j\in\Lambda_2}k'_j(\bG).
\label{eq:addup}
\ee

\paragraph{Constraints on the top and the bottom layer.}
We first fix the degree sequence on both layers, i.e., we constrain $\vec{k}_{1\to2}(\bG)$ and 
$\vec{k}_{2\to1}(\bG)$ to the values $\vec{k}_{1\to2}^*=(k^*_i)_{i\in\Lambda_1}$ and 
$\vec{k}_{2\to1}^*=(k'^*_i)_{i\in\Lambda_2}$ respectively. The constraints are therefore
\be
\label{twocon}
\vC^* =\{\vec{k}_{1\to2}^*,\vec{k}_{2\to1}^*\}.
\ee 
As mentioned before, we allow $n_1$ and $n_2$ to depend on $n$, i.e., $n_1=n_1(n)$ 
and $n_2=n_2(n)$. In order not to overburden the notation, we suppress the dependence 
on $n$ from the notation.

We abbreviate
\be
\begin{aligned}
&m^*=\max_{i\in\Lambda_1} k^*_i, \quad 
m'^*=\max_{j\in\Lambda_2} k'^*_j,\\
&f_{1\to 2}^{({n_1})}={n_1}^{-1} \sum_{i\in\Lambda_1} \delta_{k_i^*}, \quad
f_{2\to 1}^{({n_2})}={n_2}^{-1} \sum_{j\in\Lambda_2}\delta_{k'^*_j},
\end{aligned}
\ee
and assume the existence of
\be
\label{Alims}
A_1=\lim_{n\to \infty} \frac{{n_1}}{{n_1}+{n_2}}, 
\quad A_2=\lim_{n\to \infty} \frac{n_2}{{n_1}+{n_2}}.
\ee
(This assumption is to be read as follows: choose $n_1=n_1(n)$ and $n_2=n_2(n)$ in such a way 
that the limiting fractions $A_1$ and $A_2$ exist.) The \emph{sparse regime} corresponds to
\be
\label{maxdegree}
m^*m'^*=o({L^*}^{2/3}), \qquad n\to\infty.
\end{equation}
We further assume that there exist $f_{1\to 2},f_{2\to 1} \in \cP(\N_0)$ such that 
\be
\label{convergence}
\lim_{{n}\to\infty} \|f_{1\to 2}^{({n_1})}-f_{1\to 2}\|_{\ell^1(g)} = 0, \quad
\lim_{{n}\to\infty} \|f_{2\to 1}^{({n_2})}-f_{2\to 1}\|_{\ell^1(g)} = 0.
\ee
The specific relative entropy is
\be
s_{{n_1}+{n_2}} = \frac{S_{{n_1}+{n_2}}(\Pmic \mid \Pcan) }{{n_1}+{n_2}}.
\ee

\begin{theorem}
\label{result}
Subject to \eqref{twocon} and \eqref{Alims}--\eqref{convergence},
\be
s_\infty = \lim_{n\to \infty} \frac{S_{{n_1}+{n_2}}(\Pmic \mid \Pcan)}{{n_1}+{n_2}}
= A_1\,\|f_{1\to 2}\|_{\ell^1(g)} + A_2\,\|f_{2\to 1}\|_{\ell^1(g)}.
\ee
\end{theorem}

\noindent 
Since $A_1+A_2=1$, it follows that $s_{\infty}>0$, so in this case ensemble equivalence 
never holds.

\paragraph{Constraints on the top layer only.}
We now partly relax the constraints and only fix the degree sequence $\vec{k}_{1\to2}(\bG)$ 
to the value
\be
\label{topcon}
\vC^* = \vec{k}^*_{1\to 2} = \big(k_i^*\big)_{i\in\Lambda_1},
\ee 
while leaving $\vec{k}_{2\to1}(\bG)$ unspecified (apart for the condition \eqref{eq:addup}).
The microcanonical number of graphs satisfying the constraint is
\be
\Omega_{\vec{k}^*_{1\to 2}} = \prod_{i\in\Lambda_1} {{n_2} \choose k_i^*} .
\label{eq:OmegaRG}
\end{equation}
The canonical ensemble can be obtained from \eqref{eq:PC} by setting 
\begin{equation}
\label{hamiltonianbip}
H(\bG,\vt)=\vt\cdot \vec{k}_{1\to 2}(\bG).
\end{equation}
Setting $\vec{\theta}=\vec{\theta^*}$ in order that equation \eqref{eq:PC} is satisfied, we can 
write the canonical probability as 
\be
\Pcan(\bG) = \prod_{i\in\Lambda_1}({p^*_i})^{k_i(\bG)} (1-{p^*_i})^{{n_2} - k_i(\bG)}
\label{eq:PCRG}
\ee
with $p^*_i = \frac{k_i^*}{n_2}$. Let
\be
f_{n_1} = {n_2}^{-1} \sum_{i\in\Lambda_2} \delta_{k^*_i} \in \cP(\N_0).
\ee
Suppose that there exists an $f \in \cP(\N_0)$ such that
\be
\label{twofixedcon}
\lim_{{n}\to\infty} \|f_{n_1}-f\|_{\ell^1(g)} = 0.
\ee

The relative entropy per node can be written as
\be
\label{relentropy}
s_{{n_1}+{n_2}} = \frac{S_{{n_1}+{n_2}}(\Pmic \mid \Pcan) }{{n_1}+{n_2}}
=\frac{{n_1}}{{n_1}+{n_2}}\|f_{n_1}\|_{\ell^1(g_{n_2})},
\ee
with 
\be g_{n_2}(k) = -\log \left[\mathrm{Bin}\left({n_2},\tfrac{k}{n_2}\right)(k)\right]
\mathbb{I}_{0\leq k \leq {n_2}}, \qquad k \in \N_0,
\ee
and $\mathrm{Bin}({n_2},\tfrac{k}{n_2})(k) = {{n_2} \choose k} (\tfrac{k}{n_2})^{k}
(\tfrac{{n_2}-k}{k})^{{n_2}-k}$ for $k=0,\ldots,{n_2}$ and equals to $0$ for $k>{n_2}$. 
We follow the convention $0\log(0)=0$.

In this partly relaxed case, different scenarios are possible depending on the specific 
realisation of the thermodynamic limit, i.e., on how  ${n_1},{n_2}$ tend to infinity. The 
ratio between the sizes of the two layers $c = \lim_{n\to\infty} \frac{{n_2}}{{n_1}}=\frac{A_2}{A_1}$ 
plays an important role.

\begin{theorem}
\label{configurationmodel}
Subject to \eqref{topcon} and \eqref{twofixedcon}:\\
{\rm (1)} If ${n_2}\to^{n\to \infty}\infty$ with ${n_1}$ fixed ($c=\infty$), 
then $s_\infty = \lim_{{n} \to \infty} s_{{n_1}+{n_2}}=0$.\\
{\rm (2)} If ${n_1},{n_2}\to^{n\to \infty}\infty$ with $c=\infty$, then 
$s_\infty = \lim_{n \to \infty} s_{{n_1}+{n_2}}=0$.\\
{\rm (3)} If ${n_1}\to^{n\to \infty}\infty$ with ${n_2}$ fixed ($c=0$), then
\be
\label{Case 3}
s_\infty = \lim_{{n} \to \infty} s_{{n_1}+{n_2}} = \|f\|_{\ell^1(g_{n_2})}.
\ee
{\rm (4)} If ${n_1},{n_2}\to^{n\to \infty}\infty$ with $c \in [0,\infty)$, then 
\begin{equation}
\label{Case 4}
s_\infty = \frac{1}{1+c} \|f\|_{\ell^1(g)}.
\end{equation}
\end{theorem}

\paragraph{Constraints on the total number of links only.}
We now fully relax the constraints and only fix the total number of links, i.e.,
\be
\label{twototal}
\vC^* = L^*. 
\ee
In analogy with the corresponding result for the uni-partite case (Theorem \ref{biplinks}), in 
this case ensemble equivalence is restored.

\begin{theorem}
\label{resultbip}
Subject to \eqref{twototal}, the specific relative entropy equals $s_{\infty}=0$.
\end{theorem}


\subsubsection{Multiple layers}
\label{S1.3.3}

We now come to our most general setting where we fix a finite number $M \in \N$ of layers. 
Each layer $s$ has $n_s$ nodes, with $\sum_{s=1}^Mn_s=n$. Let $v_i^{(s)}$ denote the 
$i$-th node of layer $s$, and $\Lambda_s=\lbrace v^{(s)}_{1},\ldots,v_{n_s}^{(s)}\rbrace$ 
denote the set of nodes in layer $s$. We may allow links \emph{both within and across} 
layers, while constraining the numbers of links among different layers separately. But we 
may as well switch off links inside or between (some of the) layers. The actual choice can 
be specified by a superstructure, which we denote as the \emph{master graph} $\bGamma$, 
in which self-loops are allowed but multi-links are not. The nodes set of $\bGamma$ is 
$\left\lbrace 1,\dots,M\right\rbrace$ and the associated adjacency matrix has entries
\be
\gamma_{s,t}(\bGamma)=
\begin{cases}
1 \qquad \mbox{if a link between layers\ } s \mbox{\ and\ } t \mbox{ exists}
\\
0 \qquad \mbox{otherwise}.
\end{cases}
\ee
The chosen set of all multi-layer graphs with given numbers of nodes, layers, and admissible 
edges (we admit edges only between layers connected in the \emph{master graph}) is 
$\cG^\sharp_n=\cG_{n_1,\ldots,n_M}(\bGamma)\subseteq \cG_n$. In \ref{sec:empirical} we 
discuss various empirically relevant choices of $\bGamma$ explicitly, while here we keep 
our discussion entirely general.

Given a graph $\bG$, for each pair of layers $s$ and $t$ (including $s=t$) we define the 
\emph{$t$-targeted degree sequence} of layer $s$ as $\vec{k}_{s\to t}(\bG)
=\big(k_i^{t}(\bG)\big)_{i\in\Lambda_s}$, where $k^t_{i}(\bG) = \sum_{j\in\Lambda_t}
g_{i,j}(\bG)$ is the number of links connecting node $i$ to all other nodes in layer $t$. 
For each pair of layers $s$ and $t$ such that $\gamma_{s,t}(\bGamma)=1$, we enforce 
the value ${\vec{k}}^{\,*}_{s\to t} =\big(k_i^{*\,t}\big)_{i\in\Lambda_s}$ as a constraint for 
the $t$-targeted degree sequence of layer $s$. For $\gamma_{s,t}(\bGamma)=0$ we have 
${\vec{k}}^{\,*}_{s\to t} =\vec{0}$, but this constraint is automatically enforced by the master 
graph. Thus, the relevant constraints are 
\be
\label{conmulti}
\vC^*=\left\lbrace {\vec{k}}^{\,*}_{s\to t}\colon\, s,t={1,\ldots,M}
\ \ \gamma_{s,t}(\bGamma)=1 \right\rbrace.
\ee
We abbreviate
\be
\label{eq:Lst}
L_{s,t}^* =\sum_{i\in \Lambda_s} k_i^{*\,t} = \sum_{j\in  \Lambda_t} k_j^{*\,s},
\quad m_{s\to t}^*=\max_{i\in \Lambda_s} k_i^{*\,t},
\quad f_{s\to t}^{(n_s)}= n_s^{-1} \sum_{i\in \Lambda_s} \delta_{k_i^{*\,t}},
\ee
where $L_{s,t}^* $ is the number of links between layers $s$ and $t$ (note that $L_{s,s}^*$ 
is \emph{twice} the number of links inside layer $s$), and assume the existence of
\be
\label{Ahlim}
A_s=\lim_{n_1,\ldots ,n_M \to \infty}\frac{n_s}{n} \quad \forall\,s,
\ee
where $\sum_{s=1}^M A_s=1$. (As before, this assumption is to be read as follows: choose 
$n_s=n_s(n)$, $1 \leq s \leq M$, in such a way that the limiting fractions $A_s$,$1 \leq s \leq M$, 
exist.) The \emph{sparse regime} corresponds to
\be
\label{maxdegreecomm}
\begin{array}{lll}
&m_{s\to t}^*m_{t\to s}^* = o({L_{s,t}^*}^{{2}/{3}}), &n_s,n_t\to\infty\,\,\forall\,s\neq t,\\
&m_{s\to s}^* = o({n_s^{{1}/{2}}}), &n_s\to\infty\,\,\forall\,s.
\end{array}
\ee
We further assume that there exists $f_{s\to t} \in \cP(\N_0)$ such that
\be
\label{convergencecomm}
\lim_{n_s\to\infty} \|f_{s\to t}^{(n_s)}-f_{s\to t}\|_{\ell^1(g)},\qquad \lim_{n_s\to\infty} 
\|f_{s\to s}^{(n_s)}-f_{s\to s}\|_{\ell^1(g)} = 0. 
\ee

\begin{theorem}
\label{resultcomm}
Subject to \eqref{conmulti} and \eqref{Ahlim}--\eqref{convergencecomm},
\be
\label{sMG}
s_{\infty} =\sum_{\substack{s,t=1\\\gamma_{s,t}(\bGamma)=1}}^M A_s\,\| f_{s\to t}\|_{\ell^1(g)}.
\ee
\end{theorem}

\noindent 
The above result shows that, unless $A_s=0$ whenever $\gamma_{s,t}(\bGamma)=1$ (i.e., unless 
only the nodes of the master graph that have no links or self-loops contribute a finite fraction of nodes 
in the corresponding layers), ensemble equivalence does not hold.


\subsubsection{Relaxing constraints in the multilayer case}
\label{S1.3.4}

We next study the effects of relaxing constraints. This deserves a separate discussion, 
since in the multi-partite setting there are more possible ways of relaxing the constraints 
than in the uni-partite and bi-partite settings. 

\paragraph{One class of layers.}
We first fix \emph{two kinds of constraints}: (1) the total number of links between some pairs 
of layers; (2) the degree sequence between some other pairs of layers. We define the set of 
the edges of the \emph{master graph} as ${\cal{E}}=\left\lbrace (s,t)\in (M\times M)\colon\,
\gamma_{s,t}(\bGamma)=1\right\rbrace$. Then, we partition ${\cal{E}}$ into two parts, namely 
$\mathcal{D},\mathcal{L}\subseteq {\cal{E}}$, with $\mathcal D\cap \mathcal L=\emptyset$, 
$\mathcal D$ and $\mathcal L$ symmetric, by requiring that $(s,t)\in \mathcal D$ ($\in \mathcal L$) 
when $(t,s)\in\mathcal D$ ($\in \mathcal L$). For each pair of layers $(s,t)\in \mathcal D$ we fix 
the degree sequence $\vec{k}^{\,*}_{s\to t}$ of every node of $\Lambda_s$ linking to $ \Lambda_t$. 
As before, we impose that $\sum_{i\in \Lambda_s} k_i^{*\,t}=\sum_{j\in  \Lambda_t} k_j^{*\,s}$. 
For each pair of layers $(s,t)\in \mathcal L$ we fix the total number of links $L_{s,t}^*$ ($L_{s,t}^*
=L_{t,s}^*$).

The effect of relaxing some constraints affects the specific relative entropy: this will decrease 
because the pairs of layers with relaxed constraints (i.e., the pairs in $\mathcal L$) no longer 
contribute.

\begin{theorem}
\label{resultcommrelax}
Subject to the above relaxation,
\be
s_\infty=\sum_{(s,t)\in\mathcal D}
 A_s\,\|f_{s\to t}\|_{\ell^1(g)} .
\ee
\end{theorem}

\noindent
In particular, equivalence holds if and only if $\mathcal D=\emptyset$ or $A_s=0$ for all $s$ 
endpoints of elements in ${\cal{E}}$. Note that, if $\mathcal D=\emptyset$, then we have a finite 
number of constraints (at most $M^2$), and this implies equivalence of the ensembles.

\paragraph{Two classes of layers.}
We may further generalise Theorem~\ref{resultcomm} as follows. Suppose that we have 
two classes of layers, ${{\cal{M}}_1}$ and ${{\cal{M}}_2}$. For every pair of layers $s,t\in 
{{\cal{M}}_1}$ such that $\gamma_{s,t}(\bGamma)=1$, we fix the degree sequences 
${\vec{k}}^{\,*}_{s\to t}$ and $ {\vec{k}}^{\,*}_{t\to s}$. For every pair of layers $s\in {{\cal{M}}_1}$, 
$t\in {{\cal{M}}_2}$, $\gamma_{s,t}(\bGamma)=1$ we fix the degree sequence ${\vec{k}}^{\,*}_{s\to t}$ 
from the layer in ${{\cal{M}}_1}$ to the layer in ${{\cal{M}}_2}$ (but not vice versa). We show 
that the resulting specific relative entropy is a mixture of the one in Theorem~\ref{resultcomm} 
and the one in Theorem~\ref{configurationmodel}. For $s=1,\dots,M$ we set $A_s =
\lim_{n_1,n_2,\dots,n_{M}\to\infty} \frac{n_s}{n}$.

\begin{theorem}
\label{resultcommgen}
Subject to the above relaxation,
\be
\begin{aligned}
s_\infty & =\sum_{\substack{s\in {{\cal{M}}_1},\  
t\in {{\cal{M}}_1\cup{\cal{M}}_2}\\ \gamma_{s,t}(\bGamma)=1}}
 A_s\,\|f_{s\to t}\|_{\ell^1(g)}.
\end{aligned}
\ee
\end{theorem}

\noindent
In particular, 
\be
s_\infty=0 \quad \Longleftrightarrow \quad A_s=0\,\,\forall\,s\in \big\lbrace u\in {\cal{M}}_1\colon\, 
\exists\, t\in {\cal{M}}_1\cup{\cal{M}}_2 \ \mbox{with\ } \gamma_{u,t}(\bGamma)=1\big\rbrace.
\ee

\paragraph{Another way for relaxing constraints.}
We may think about another way for relaxing the constraints.  We assume that $\gamma_{s,t}(\bGamma)
=1$ for all $s,t=1,2,\dots,M$ and we fix ${\vec{k}}^{\,*}_{s}=\sum_{t=1}^M{\vec{k}}^{\,*}_{s\to t}$ for 
each $s=1,2,\dots,M$. This means that for each node we fix its degree sequence (no matter to which 
target layer, possibly its own layer). In this case we lose the multi-layer structure: constraints are no 
longer involving pairs of layers and the graphs are effectively uni-partite. This is the same case described 
in the configuration model of Theorem~\ref{resultunip}. There are still an extensive number of local 
constraints, and the ensembles are nonequivalent.


\section{Discussion}
\label{S2}

In this section we discuss various important implications of our results. We first consider properties 
that are fully general, and afterwards focus on several special cases of empirical relevance.


\subsection{General considerations}
\label{sec:general}


\paragraph{Poissonisation.}
The function $g$ in \eqref{gdef} has an interesting interpretation, namely,
\be
g(k) = S\big(\,\delta[k] \mid \mathrm{Poisson}[k]\,\big)
\ee
is the relative entropy of the Poisson distribution with average $k$ w.r.t.\ the Dirac distribution
with average $k$. The specific relative entropy in \eqref{resultunip} for the uni-partite setting
can therefore be seen as a sum over $k$ of contributions coming from the nodes with 
fixed, respectively, average degree $k$. The microcanonical ensemble forces the degree 
of these nodes to be exactly $k$ (which corresponds to $\delta[k]$), while the canonical 
ensemble, under the sparseness condition in \eqref{eq:sparse1}, forces their degree to be 
Poisson distributed with average $k$. The same condition ensures that in the limit as 
$n\to\infty$ the constraints act on the nodes essentially independently.  

The same interpretation applies to Theorems~\ref{result}--\ref{configurationmodel} and
\ref{resultcomm}--\ref{resultcommgen}. The result in Theorem~\ref{configurationmodel}(3)
shows that in the bi-partite setting, when one of the layers tends to infinity while the other 
layer does not, Poissonisation does not set in fully. Namely, we have
\be
s_n =  \sum_{k=1}^n f(k) g_n(k), \qquad g_n(k) 
= S\big(\,\delta[k] \mid \mathrm{Bin}(n,\tfrac{k}{n})\,\big).
\ee
In words, the canonical ensemble forces the nodes in the infinite layer with average degree 
$k$ to draw their degrees towards the $n$ nodes in the finite layer essentially independently, 
giving rise to a binomial distribution. Only in the limit as $n\to\infty$ does this distribution 
converge to the Poisson distribution with average $k$. 

\paragraph{Additivity vs.\ non-additivity.}
In all the other examples known so far in the literature, the generally accepted explanation 
for the breaking of ensemble equivalence is the presence of a non-additive energy, induced 
e.g.\ by long-range interactions \cite{Ruffo}, \cite{Ruffo2}. However, in the examples 
considered in the present paper, nonequivalence has a different origin, namely, the presence 
of an extensive number of local constraints. As we now show, this mechanism is completely 
unrelated to non-additivity and is therefore a novel mechanism for ensemble nonequivalence.

Intuitively, the energy of a system is additive when, upon partitioning the units of the 
system into non-overlapping subunits, the `interaction' energy between these subunits is 
negligible with respect to the internal energy of the subunits themselves. The `physical' 
size of the systems considered in this paper is given by the number $n$ of nodes, i.e., we 
are defining the network to become `twice as large' when the number of nodes is doubled. 
Think, for instance, of a population of $n$ individuals and the corresponding social network 
connecting these individuals: we say that the size of the network doubles when the population 
doubles. Consistently, in \eqref{eq:sn} we have defined the specific relative entropy $s_n$ 
by diving $S_n$ by $n$. In accordance with this reasoning, in order to establish whether in 
our systems ensemble equivalence has anything to do with energy additivity, we need to 
define the latter \emph{node-wise}, i.e., with respect to partitioning the set of nodes into 
nonoverlapping subsets. Note that, in the presence of more than one layer, we have allowed 
for the number of nodes in some layer(s) to remain finite (in general, to grow subextensively) 
as the total number of nodes goes to infinity (see for instance Theorem~\ref{configurationmodel}). 
In such a situation it makes sense to study additivity only with respect to the nodes in those 
layers that are allowed to grow extensively in the thermodynamic limit. 

Formally, if we let $\mathcal{I}$ denote the union of all layers for which $A_s>0$ (see~\eqref{Ahlim}), 
then we say that the energy is \emph{node-additive} if the Hamiltonian~\eqref{eq:H} can be written 
as 
\begin{equation}
H(\bG, \vt) = \sum_{i\in\mathcal{I}}H_i(\bG,\vt)\qquad\forall\, \bG\in\cG^\sharp_n,
\label{eq:additivity}
\end{equation}
where the $\{H_i\}_{i\in\mathcal{I}}$ do not depend on common subgraphs of $\bG$ (i.e., each 
of them can be restricted to a distinct subgraph of $\bG$), and are therefore independent random 
variables.

The case of uni-partite graphs with fixed degree sequence (Theorem~\ref{resultunip}) is an 
example of ensemble \emph{nonequivalence} with \emph{non-additive} Hamiltonian, because 
the latter is defined as $H(\bG, \vt) = \sum_{i=1}^n \theta_i k_i(\bG)$ and cannot be rewritten in 
the form of \eqref{eq:additivity} with independent $\{H_i(\bG,\vt)\}$: the degrees $k_i(\bG)$ and 
$k_j(\bG)$ of any two distinct nodes $i$ and $j$ depend on a common subgraph of $\bG$, i.e., 
the dyad $g_{i,j}(\bG)$. In the example of uni-partite graphs with a fixed total number of links 
(see \eqref{linkcon}), the energy has the form $H(\bG, \vt) = \theta L(\bG)= \tfrac12\theta 
\sum_{i=1}^n k_i(\bG)$, which is still \emph{non-additive}. However, the ensembles are in this case 
\emph{equivalent} (see Theorem~\ref{biplinks}).

By contrast, the case of bi-partite graphs with fixed degree sequence on the top layer and the 
nodes in the other layer growing subextensively (case (3) of Theorem~\ref{configurationmodel}) 
is an example of ensemble \emph{nonequivalence} with an \emph{additive} Hamiltonian. Indeed, 
from \eqref{hamiltonianbip} we see that $H(\bG,\vt)$ is now a linear combination of the $n_1$ 
degrees of the nodes in layer $\Lambda_1$, each of which depends only on the (bi-partite) 
subgraph obtained from the corresponding node of the top layer and all the nodes of the bottom 
layer. Here, unlike the uni-partite case, all these subgraphs are disjoint. Despite being node-additive, 
when $A_1=1$ ($c=0$) this Hamiltonian leads to nonequivalence, as established in~\eqref{Case 3}. 
Similar examples can be engineered using some of the relaxations in Section~\ref{S1.3.4}. Finally, 
the case of bi-partite graphs with fixed total number of links (Theorem~\ref{resultbip}) is an 
example of ensemble \emph{equivalence} with an \emph{additive} Hamiltonian.

The four examples above show that additivity or non-additivity of the Hamiltonian does \emph{not} 
influence the breaking of ensembles equivalence in the examples considered here. What 
matters is the \emph{extensiveness} of the number of constraints. This observation was 
already made in \cite{SdMdHG15}, and is confirmed in full generality for the multi-layer 
setting treated in the present paper. Indeed, our results indicate that, whenever the number 
$\kappa$ of constraints on the degrees is \emph{subextensive}, i.e., $\kappa=o(n)$ where 
$n$ is the number of nodes, ensemble equivalence is restored.

Note that the above notion of \emph{node additivity} should not be confused with that of 
\emph{edge additivity}, i.e., the fact that the Hamiltonian can be written as a sum over independent 
pairs of nodes. Due to the linearity of the chosen (local) constraints on the entries $\{g_{i,j}\}_{i,j=1}^n$ 
of the adjacency matrix of the graph $\bG$, our examples are always edge-additive (irrespective 
of whether they are ensemble-equivalent), while they may or may not be node-equivalent, as we 
have seen. In either case, there is no relation between additivity and equivalence.

We stress again that the extensivity of the (local) constraints is, with respect to the mechanisms 
for nonequivalence already explored in the literature so far, an additional (and previously 
unrecognised) \emph{sufficient} mechanism. It is obviously not the only one, and definitely 
\emph{not a necessary one}, as exemplified by the fact that, in dense networks, nonequivalence 
has been found even in the presence of only two constraints, such as the total numbers of edges 
and triangles~\cite{RS,RS3,RS2,chatterjee}. However, while in the previous examples the breaking 
of equivalence arises from the nonlinearity (with respect to $\{g_{i,j}\}$) of some constraint and is 
typically found in a specific (usually critical) region of the parameter space separating phases
where ensemble equivalence still applies, in our setting ensemble nonequivalence arises from 
the extensiveness of the number of (linear) constraints and extends to the entire space of parameters 
of the models. In this sense it is a stronger form of nonequivalence. Moreover, while the nonequivalence 
of network ensembles with a finite number of constraints was previously reported only for dense 
graphs, we are documenting it for the unexplored regime of sparse graphs. 

\paragraph{A principled choice of ensembles.}
Ensembles of random graphs with constraints are used for many practical purposes. Two 
important examples are \emph{pattern detection} and \emph{network reconstruction}. For 
concreteness, we briefly illustrate these examples before we emphasize the implications 
that our results have for these and other applications.

Pattern detection is the identification of nontrivial structural properties in a real-world network, 
through the comparison of such network with a suitable null model \cite{mymethod}. For instance, 
\emph{community detection} is the identification of groups of nodes that are more densely 
connected with each other than expected under a null model \cite{community_santo}, 
\cite{community_porter}
 (in 
Section~\ref{sec:empirical} we discuss the relation between our models and community detection 
in more detail). A null model is a random graph model that preserves some simple topological
properties of the real network (typically local, like the degree sequence) and is otherwise 
completely random. So, maximum-entropy ensembles of graphs with given degrees are a 
key tool for pattern detection. 

Network reconstruction employs purely local topological information to infer the higher-order 
structural properties of a real-world network \cite{myenhanced}. This problem arises whenever 
the complete structure of a network is not known (for instance, due to confidentiality or privacy 
issues), but local properties are. An example relevant for the epidemiology of sexually transmitted 
diseases is the network of sexual contacts among people, for which only aggregate information 
(the total number of contacts with different partners) can be typically surveyed in a population.
In such cases, optimal inference about the network can be achieved by maximising the entropy 
subject to the known (local) constraints, which again leads to the ensembles with fixed degrees 
considered here. 

The aforementioned applications, along with similar ones, make use of random graphs with local 
constraints. Our proof of nonequivalence of the corresponding ensembles have the following 
important implication. While for ensemble-equivalent models it makes practically no difference 
whether a microcanonical or canonical implementation is applied to large networks, for nonequivalent 
models different choices of the ensemble lead to asymptotically different results. As a consequence, 
while for applications based on ensemble-equivalent models the choice of the working ensemble 
can be arbitrary or be done on mathematical convenience (as usually done), for those based on 
nonequivalent models the choice should be principled, i.e., dictated by a theoretical criterion that 
indicates \emph{a priori} which ensemble is the appropriate one.

Among the possible criteria, we suggest one that we believe appropriate whenever the available 
data are subject to (even small) errors, i.e., when the measured value $\vec{C}^*$ entering as 
input in the construction of the random graph ensemble is, strictly speaking, the best available 
estimate for some unknown `true' (error-free) value $\vec{C}^{\times}$. In this situation, we want 
that possible small deviations of $\vec{C}^{*}$ from $\vec{C}^{\times}$ result in small devations 
of $P^*_\textrm{mic}$ and $P^*_\textrm{can}$ from the corresponding $P^\times_\textrm{mic}$ 
and $P^\times_\textrm{can}$. Now, if $\vec{C}^{*}\ne\vec{C}^{\times}$ (no matter how ``small'' 
and in which norm this difference is taken), then $P^*_\textrm{mic}$ will attach zero probability 
to any graph $\bG^\times$ that realises the `true' constraint $\vec{C}^\times$: $P^*_\textrm{mic}
(\bG^\times)=0$, while $P^\times_\textrm{mic}(\bG^\times)\ne 0$. Indeed, $P^*_\textrm{mic}$ 
and $P^\times_\textrm{mic}$ will have non-overlapping supports, so they will sample distinct 
sets of graphs. This means that even small initial errors in the knowledge of the constraints will 
be severely propagated to the entire microcanonical ensemble, and inference based on the latter 
will be highly biased. In particular, the `true' network will never be sampled by $P^*_\textrm{mic}$.
On the other hand, if the difference between $\vec{C}^{*}$ and $\vec{C}^{\times}$ is small, then 
the difference between $P^*_\textrm{can}$ and $P^\times_\textrm{can}$ will also be small. So, 
even though $\vec{C}^{\times}$ is unknown, any graph $\bG^\times$ that realises this value will 
be given a probability $P^*_\textrm{can}(\bG^\times)$ that is nonzero and not very different from 
the probability $P^\times_\textrm{can}(\bG^\times)$ that would be obtained by knowing the true 
value $\vec{C}^{\times}$. In general, small deviations of $\vec{C}^{*}$ from $\vec{C}^{\times}$ 
imply that $P^*_\textrm{can}(\bG)$ is not very different from $P^\times_\textrm{can}(\bG)$ for 
any graph $\bG$, as desired. This implies that \emph{even if $\vec{C}^{*}$ is affected by small 
errors, then a principled choice of ensembles is the canonical one.} So, besides being the 
mathematically simpler option, we argue that canonical ensembles are also the most appropriate 
choice in the presence of `noise'. A similar claim was already made in \cite{myunbiased}, and is 
here strengthened by our proof of nonequivalence.


\subsection{Special cases of empirical relevance}
\label{sec:empirical}

Different choices of the master graph $\bGamma$ induce different structural features in the graphs 
of the ensemble $\cG^\sharp_n$. Convenient choices allow us to consider certain classes of graphs 
that have been introduced recently to study appropriate types of real-world networks of empirical 
relevance. We discuss some of these choices below. The full generality of our results in 
Section~\ref{S1.3.3} allows us to immediately draw conclusions about the (non)equivalence of the 
corresponding ensembles in each case of interest. As an important outcome of this discussion, all 
the empirically relevant ensembles of graphs turn out to be nonequivalent. In line with our general 
observation at the end of the previous section, this implies that a principled choice of ensembles is 
needed in all practical applications.

\paragraph{Scale-free uni-partite networks.}
Clearly, the trivial case when the master graph has a single node ($M=1$) with a self-loop, i.e., 
$\gamma_{11}(\bGamma)=1$, corresponds to the class of uni-partite graphs we considered in 
Section~\ref{S1.3.1}. Many real-world networks, at least at a certain level of aggregation, admit 
such uni-partite representation. Examples include the Internet, the World Wide Web and many 
biological, social and economic networks. A common property displayed by most of these 
real-world networks is the presence of a ``broad'' empirical degree distribution, often consistent 
with a power-law distribution with an upper cut-off \cite{cutoff}. Networks with a power-law degree 
distribution are said to be \emph{scale-free} \cite{guidosbook}. This empirical observation implies 
that real-world networks are very different from Erd\H{o}s-R\'enyi random graphs (which have a 
much narrower degree distribution) and are more closely reproduced by a configuration model 
with a truncated power-law degree distribution $f_n$ (see \eqref{empdef}) of the form $f_n(k) 
= A_{\gamma,n} k^{-\gamma}\mathbb{I}_{1 \leq k \leq k_c(n)}$ with $\gamma>1$, $A_{\gamma,n}$ 
the normalisation constant, and $\lim_{n\to\infty} k_c(n)=\infty$ and $k_c(n) = o(\sqrt{n})$. The 
so-called \emph{structural cut-off} $k_c(n)$ makes the networks sparse, as in condition 
\eqref{eq:sparse1} \cite{cutoff}. Since $\lim_{n\to\infty} \|f_n-f\|_{\ell^1(g)}=0$ with $f(k) 
= k^{-\gamma}/\zeta(\gamma)$ for $k\ge 1$ and $0$ elsewhere, where $\zeta$ is the 
Riemann zeta-function, our result in \eqref{sec:unis} tells us that
\begin{equation}
s_\infty = \sum_{k\in\N} g(k)\,f(k) =  \frac{1}{\zeta(\gamma)} \sum_{k\in\N} g(k)\,k^{-\gamma}.
\end{equation}
Since $g(k) = \tfrac12\log (2\pi k)+O(k^{-1})$ as $k\to\infty$, we find that $s_\infty$ tends to 1 
as $\gamma \to\infty$ and diverges like $\sim 1/2(\gamma-1)$ as $\gamma \downarrow 1$. 
This result shows that the simplest random graph ensemble consistent with the scale-free 
character of real-world networks is nonequivalent. Interestingly, as the tail exponent $\gamma$ 
decreases, the degree distribution becomes broader and the degree of nonequivalence 
increases. A similar conclusion was drawn in \cite{SdMdHG15}.

\begin{remark}
{\rm Suppose that for each $n\in\N$ the degrees are drawn in an i.i.d.\ manner from the 
truncated degree distribution $f_n$. Suppose further that $\sum_{k\in\N_0} kf(k)<\infty$,
i.e., $\gamma>2$. Then, because $\sup_{n\in\N} \sum_{k\in\N_0} kf_n(k) = \sum_{k\in\N_0} 
kf(k)<\infty$, conditional on the sum of the degrees being even, the degree sequence is 
graphical with a probability tending to one as $n\to\infty$. This fact is the analogue of the 
result in \cite{arratialiggett2005} mentioned in Remark~\ref{remAL}, and its proof is a 
straightforward extension of the argument in \cite{arratialiggett2005}. Truncation improves 
the chance of being graphical.}
 \end{remark}

\paragraph{Multipartite networks.}
The case when the master graph has only $M=2$ interconnected nodes \emph{and no self-loops}, 
i.e., $\gamma_{1,2}(\bGamma)=\gamma_{2,1}(\bGamma)=1$ and $\gamma_{1,1}(\bGamma)
=\gamma_{2,2}(\bGamma)=0$, coincides with the class of bi-partite graphs discussed in 
Section~\ref{S1.3.2}. Popular real-world examples relevant to economics, ecology and scientometrics 
are bank-firm, plant-pollinator and author-paper networks, respectively. In this case as well, 
empirical evidence shows that real-world bi-partite networks have broad degree distributions (at 
least on one of the two layers, and typically on both). Random graph models with only a global 
constraint on the total number of links (as in Theorem~\ref{resultbip}) are therefore unrealistic. 
The minimal ensemble that is consistent with the properties of most real-world bi-partite networks 
requires the specification of the degree sequence(s) as constraint(s) and is therefore nonequivalent.

A direct generalisation of the bi-partite case is when $\bGamma$ is an $M$-dimensional matrix 
with zeroes along the diagonal and ones off the diagonal: $\gamma_{s,s}(\bGamma)=0$ $\forall s$ 
and $\gamma_{s,t}(\bGamma)=1$ for all $s \ne t$. The induced graphs in $\cG^\sharp_n$ have 
an ``all-to-all'' multipartite structure (i.e., links are allowed between all pairs of distinct layers, but 
not inside layers). From our Theorem~\ref{resultcomm} it follows that if the $t$-targeted degree 
sequences are specified as a constraint, then the relative entropy in the all-to-all multipartite case 
is
\be
\label{alltoall}
s_\infty =\sum_{\substack{s,t=1\\s\ne t}}^M A_s\,\| f_{s\to t}\|_{\ell^1(g)}>0,
\ee
which proves again ensemble nonequivalence.

\paragraph{Stochastic block-models.}
Another important example is when the master graph is a complete graph with all self-loops 
realised, i.e., $\gamma_{s,t}(\bGamma)=1$ for all $s,t$. This prescription generates the class 
of so-called \emph{stochastic block-models}, which are very popular in the social network 
analysis literature \cite{simpleBM}, \cite{degreecorrectedBM}, \cite{fronczakBM}. The earliest 
and simplest stochastic block-model \cite{simpleBM} is one where only the total numbers of links 
between all pairs of blocks (including within each block) are specified. When we identify blocks with 
layers, this model coincides with our relaxed model considered in Theorem~\ref{resultcommrelax}, 
with $\mathcal D=\emptyset$. It follows as a corollary that this model is ensemble equivalent:
\be
s_\infty =0.
\ee
However, this model predicts that, within each block, the expected topological properties of the 
network are those of an Erd\H{o}s-R\'enyi random graph, a property that is contradicted by 
empirical evidence. So, unless the number of blocks is chosen to be comparable with the number 
of nodes (which in our case is contradicted by the requirement that $M$ is finite), the traditional 
block-model is not a good model of real-world networks.

More recently, emphasis has been put on the more realistic \emph{degree-corrected stochastic 
block-model} \cite{degreecorrectedBM}, where an additional constraint is put on the degree of all 
nodes. An even more constrained variant of this model has been proposed in \cite{fronczakBM}, 
where the constraints coincide with the $t$-targeted degree sequences $\{\vec{k}_{s\to t}\}_{s,t}$ 
among all pairs of blocks. To distinguish this model from the ``generic'' degree-corrected block-model, 
we call it the \emph{targeted degree-corrected block-model}. This coincides with our model 
in Section~\ref{S1.3.3}, with the block structure given by the (complete) master graph. From 
Theorem~\ref{resultcomm} we calculate the relative entropy as
\be
\label{degreecorrected}
s_\infty =\sum_{s,t=1}^M A_s\,\| f_{s\to t}\|_{\ell^1(g)}>0.
\ee

We can therefore conclude that, unlike the traditional block-model considered above, the targeted 
degree-corrected model is ensemble nonequivalent. We also note that, unlike stated in \cite{fronczakBM}, 
the targeted degree-corrected block-model is not just a reparametrisation of the untargeted 
degree-corrected model. While fixing the targeted degree sequences automatically realises the 
constraints of the untargeted model, the converse is not true. Being a relaxation of the targeted model, 
we expect the untargeted model to have a relative entropy smaller than in \eqref{degreecorrected}, 
further illustrating the difference between the two models. Yet, we expect the relative entropy in the 
untargeted model to be strictly positive for, every choice of the degree sequence, since there is still 
an extensive number of active constraints. This would support the claim made in \cite{entropyBM} 
that, \emph{for small values of the degrees}, the degree-corrected block-models with soft and hard 
constraints are not equivalent in the thermodynamic limit. At the same time, it would contradict the 
claim made in the same reference that, \emph{if all degrees become large} (but still in the sparse 
regime), the two ensembles become equivalent. Indeed, from the behaviour of $g(k)$ for large $k$ 
(see~\eqref{g(k)}) and the normalisation by $n$ in~\eqref{eq:sn}, we expect a finite specific relative 
entropy in that case as well.

\paragraph{Networks with community structure.}
Another very important class of graphs that are studied intensively in the literature are graphs with 
community structure \cite{community_santo}, \cite{community_porter}. This class is related to the 
block-models described above, but is in general different. Community structure is loosely defined 
as the presence of groups of nodes that are more densely interconnected internally than with each 
other. One of the possible ways to quantitatively define the presence of communities in a real-world 
network is in terms of a positive difference between the realised number of intra-community links 
and the corresponding expected number calculated under a certain null hypothesis. This definition 
can be made more explicit by introducing the concept of \emph{modularity} \cite{community_santo}, 
\cite{community_porter}. For a graph with $n$ nodes, a non-overlapping partition of nodes into $M$ 
communities can be specified by the $n$-dimensional vector $\vec{\sigma}$, where the $i$-th entry
$\sigma_i\in \{1, \dots, M\}$ is an integer number labelling the community to which node $i$ is 
assigned by that particular partition. For a given real-world graph $\bG^*$, the modularity is a 
function on the space of possible partitions, defined as
\be
\label{eq:modularity}
Q_{\bG^*}(\vec{\sigma})
= K_{\bG^*} \sum_{1\le i<j\le n} \left(g_{ij}(\bG^*)-\langle g_{ij}\rangle\right)
\delta_{\sigma_i,\sigma_j},
\ee
where $K_{\bG^*}$ is an (inessential) normalisation constant (independent of the partition $\vec{\sigma}$) 
intended to have the property $Q_{\bG^*} \in [-1,+1]$, and $\langle g_{ij}\rangle$ is the expected value 
of $g_{ij}(\bG)$ under the null hypothesis. The null hypothesis leads to a \emph{null model} for the 
real-world network ${\bG^*}$. The most popular choice for this null model is the canonical configuration 
model in the sparse regime, which gives $\langle g_{ij}\rangle=k^*_i k^*_j/2L^*$ for $i\ne j$ and 
$\langle g_{ii}\rangle=0$, where $k^*_i$, $k^*_j$ and $L^*$ are all calculated on $\bG^*$ (see 
\eqref{peq:pijCL} in the proof of Theorem~\ref{resultunip}).

Now, if the real-world network $\bG^*$ is indeed composed of communities, then the partition 
$\vec{\sigma}^\dagger$ that encodes these communities will be such that $Q_{\bG^*}(\vec{\sigma}^\dagger)
>0$, i.e., the total number of links inside communities will be larger than the expected number 
under the null model. More stringently, the `optimal' partition into communities can be defined as 
the one that maximises $Q_{\bG^*}(\vec{\sigma})$, provided that the corresponding value 
$\max_{\vec{\sigma}}{Q_{\bG^*}(\vec{\sigma})}$ is positive. Indeed, one of the most popular ways 
in which communities are looked for in real-world networks is through the process of modularity 
maximisation. The higher the value of the maximised modularity, the sharper the community structure. 
In practice, the problem of community detection is complicated by the possible existence of many 
local minima of $Q_{\bG^*}(\vec{\sigma})$ and by the fact that $Q_{\bG^*}(\vec{\sigma}^\dagger)$ 
may be positive even for ``noisy communities'', i.e., communities induced by chance only out of 
randomness in the data.

In our setting, community structure can be easily induced in the multilayer graph ensemble 
$\cG^\sharp_n=\cG_{n_1,\ldots,n_M}(\bGamma)$ through a convenient choice of the master 
graph $\bGamma$ and of the constrained $t$-targeted degree sequences $\{\vec{k}_{s\to t}^*\}$. 
First, we identify the $M$ layers $\{\Lambda_s\}$ with the desired communities and define the 
corresponding partition $\vec{\sigma}^\dagger$ through $\sigma^\dagger_i=\Lambda_s$ if 
$i\in\Lambda_s$. Next, we require that the master graph $\bGamma$ has all possible self-loops, 
plus a desired number of additional edges that need not be maximal (pairs of distinct 
communities are not necessarily connected in real-world networks). Finally, we need to require 
that the $t$-targeted degree sequences induce an excess of intra-community links with respect 
to the null model, so that the modularity is at least positive, i.e., $Q_{\bG^*}(\vec{\sigma}^\dagger)
>0$, and at best maximised by the desired partition, i.e., $\vec{\sigma}^\dagger
=\textrm{argmax}_{\vec{\sigma}}\,{Q_{\bG^*}(\vec{\sigma})}$. To this end, we rewrite
\begin{equation}
\begin{aligned}
Q_{\bG^*}(\vec{\sigma}^\dagger)
&= K_{\bG^*}\sum_{1\le i<j\le n}\left(g_{ij}(\bG^*)-\langle g_{ij}\rangle\right)
\delta_{\sigma_i^\dagger,\sigma_j^\dagger}\\
&= \frac{K_{\bG^*}}{2}\sum_{1\le i,j\le n}\left(g_{ij}(\bG^*)-\langle g_{ij}\rangle\right)
\delta_{\sigma_i^\dagger,\sigma_j^\dagger}\\
&= \frac{K_{\bG^*}}{2}\sum_{s=1}^M\sum_{i,j\in\Lambda_s}
\left(g_{ij}(\bG^*)-\frac{k_i^*k_j^*}{2L^*}\right)\\
&= \frac{K_{\bG^*}}{2}\sum_{s=1}^M\left(L^*_{s,s}-\frac{1}{2L^*}
\Big(\sum_{i\in\Lambda_s}k_i^*\Big)^2\right)\\
&= \frac{K_{\bG^*}}{2}\sum_{s=1}^M\left(L^*_{s,s}-\frac{1}{\sum_{s,t=1}^M L^*_{s,t}}
\Big(\sum_{t=1}^M L_{s,t}^*\Big)^2\right),
\label{lastline}
\end{aligned}
\end{equation}
where we use $g_{ii}(\bG^*)=\langle g_{ii}\rangle=0$, $k_i^*=\sum_{t=1}^Mk_i^{*t}$ and 
$2L^*=\sum_{s,t=1}^M L^*_{s,t}$. So, the weaker condition $Q_{\bG^*}(\vec{\sigma}^\dagger)>0$ 
is realised by requiring that $\{\vec{k}^*_{s\to t}\}$ satisfies the inequality 
\be
\label{commcondition}
\sum_{s=1}^M L^*_{s,s}
>\frac{\sum_{s=1}^M\Big(\sum_{t=1}^M L_{s,t}^*\Big)^2}{\sum_{s,t=1}^M L^*_{s,t}},
\ee
where $L^*_{s,t}=\sum_{i\in \Lambda_s}k_i^{*t}$. The above inequality explicitly states that the 
number of realised intra-community edges counted in the left-hand side should be larger than the 
expected number calculated in the right-hand side. The stronger condition $\vec{\sigma}^\dagger
=\textrm{argmax}_{\vec{\sigma}}\,{Q_{\bG^*}(\vec{\sigma})}$ should instead be enforced by 
looking for the specific $\{\vec{k}^*_{s\to t}\}$ that maximises \eqref{lastline}.

Independently of how communities are induced in our framework, our results show that 
\emph{ensembles of random graphs with community structure} (according to the definition 
above) \emph{are nonequivalent}, with a relative entropy given by \eqref{sMG} where the 
degree distributions $\{f_{s\to t}\}$ are induced by suitable $t$-targeted degree sequences 
that realise \eqref{commcondition} and possibly also $\vec{\sigma}^\dagger
=\textrm{argmax}_{\vec{\sigma}}\,{Q_{\bG^*}(\vec{\sigma})}$.

\paragraph{Multiplex networks and time-varying graphs.}
Two other important classes of graphs that have recently gained attention are those of \emph{multiplex} 
networks \cite{multiplex} and \emph{time-varying graphs} \cite{timevarying}. 

Multiplex networks are networks where the same set of nodes can be connected by $M$ different 
types of links \cite{multiplex}. Two examples, both studies in~\cite{mymultiplexity}, are the multiplex 
of international trade in different products (where nodes are world countries and links of different 
type represent international trade in different products) and the multiplex of flights by different airlines 
(where nodes are airports and links of different type represent flights operated by different companies). 
An equivalent and widely used representation for a multiplex is one where a number $M$ of 
layers is introduced, the same nodes are replicated in each layer, and inside each layer an ordinary 
graph is constructed, specified by all links of a single type. Links only exist within layers, and not 
across layers. Indeed, what `couples' the different layers and makes a real-world multiplex different 
from a collection of independent layers is the empirical fact that the topological properties of the 
layer-specific networks are typically strongly (either positively or negatively) correlated. For instance, 
networks of trade in different products have a similar structure, and most notably countries that 
are `hubs' in one layer are likely to be hubs in other layers as well. By contrast, airports that are 
hubs for a domestic airline are likely not to be hubs for other domestic airlines \cite{mymultiplexity}. 
This means that, for each node $i$ in real-world networks, the $M$ numbers of intra-layer links 
(i.e., the \emph{intra-layer degrees}) are in general (anti)correlated. 

Time-varying graphs are collections of temporal snapshots of the same network \cite{timevarying}. 
If the set of nodes in the network does not change with time, then a time-varying graph can be 
represented as a multiplex where each temporal snapshot is a single layer. (Note that multiplex 
networks themselves can vary over time \cite{gin_growmultiplex}.) Again, while not interacting 
directly via links, the different layers are mutually dependent because of empirical correlations 
between the properties of the same physical network across its temporal snapshots. Therefore 
this type of time-varying graphs can be treated in a way formally similar to that used for multiplex 
networks, the only difference being that a natural temporal ordering can be defined for the 
snapshots of time-varying graphs, while this is in general not true for the layers of a multiplex.

In our framework, a multiplex or time-varying network can be introduced by identifying each 
link type with a layer $\Lambda_s$ and by requiring that the only edges of the master graph 
$\bGamma$ are self-loops, i.e., $\gamma_{s,s}(\bGamma)=1$ for $s=1,M$ and $\gamma_{s,t}
(\bGamma)=0$ for $s\ne t$. Note that this specification, which implies $\vec{k}_{s\to t}^*=\vec{0}$ 
for $s\ne t$, is somehow `dual' to the one defining all-to-all multipartite networks (see above). 
The fact that nodes in different layers are replicas of the same set of $n$ nodes implies that 
$|\Lambda_s|$ is the same for all $s$, i.e., $n_s=n/M$. Finally, the `coupling' between the 
topological properties of different layers can be introduced by assigning (anti)correlated 
$t$-targeted degree sequences, i.e., by choosing (anti)correlated entries for every pair of 
vectors $\vec{k}_{s\to s}^*$ and $\vec{k}_{t\to t}^*$, $s\ne t$. Real-world multiplexes, including 
the two examples made above, are well reproduced by such a model \cite{mymultiplexity}. 
Our results imply that the relevant ensembles are nonequivalent. In particular, as a corollary 
of Theorem~\ref{resultcomm} we have
\be
\label{eq:multiplex}
s_\infty =\frac{1}{M}\sum_{s=1}^M \,\| f_{s\to s}\|_{\ell^1(g)}.
\ee
So, the relative entropy between the microcanonical and canonical distributions is the average of 
the relative entropy of all layers, where for each layer $s$ the relative entropy is the same as that 
obtained for a uni-partite network with $n/M$ nodes and limiting degree distribution $f_{s\to s}$ 
(see Theorem~\ref{resultunip}). Moreover, the presence of correlations between $\vec{k}_{s\to s}^*$ 
and $\vec{k}_{t\to t}^*$ translate into dependencies between $\| f_{s\to s}\|_{\ell^1(g)}$ and 
$\| f_{t\to t}\|_{\ell^1(g)}$. In particular, in case of perfect correlation ($\vec{k}_{s\to s}^*=\vec{k}_{t\to t}^*$ 
for all $s,t$), all the degree distributions are equal to a common one $f_{s\to s}=f$ $\forall s$, and 
we get
\be
\label{eq:multiplex2}
s_\infty = \,\| f\|_{\ell^1(g)}.
\ee
In this case, the degree of nonequivalence is the same as that obtained for a single uni-partite network 
with $n/M$ nodes and limiting degree distribution $f$ (see Theorem~\ref{resultunip}).

\paragraph{Interdependent multilayer networks.}
Finally, we discuss the class of \emph{interdependent multilayer networks}, which are multiplex networks 
with the addition of inter-layer links \cite{multiplex}. Nodes in different layers are still replicas of the same set 
of nodes, so we still have $n_s=n/M$ for all $s$. Similarly, the topological properties of different intra-layer 
networks are still (anti)correlated, which can be again realised by choosing (anti)correlated entries for 
every pair of vectors $\vec{k}_{s\to s}^*$ and $\vec{k}_{t\to t}^*$, $s\ne t$. However, while we still require 
$\gamma_{s,s}(\bGamma)=1$ for $s=1,M$, now we no longer require $\gamma_{s,t}(\bGamma)=0$ for 
$s\ne t$. Therefore the degree of nonequivalence can only increase with respect to \eqref{eq:multiplex}.
Indeed, Theorem~\ref{resultcomm} now leads to
\be
s_\infty=\frac{1}{M}\sum_{\substack{s,t=1\\\gamma_{s,t}(\bGamma)=1}}^M \,\| f_{s\to t}\|_{\ell^1(g)},
\ee
which shows that the relative entropy is no longer only an average over the layer-specific relative entropies, 
since inter-layer relative entropies give additional contributions.

\paragraph{Networks of networks.}
A final class of graphs worth mentioning is the so-called \emph{networks of networks}, sometimes 
constructed by different `micro-networks' that are coupled together into a `macro-network' where 
each node is a micro-network itself \cite{netonets}. This class is similar to the interdependent 
multilayer networks considered above, but here there is no identification of the nodes in different 
layers to the same physical entity. An example is provided by multi-scale transport networks, 
where different cities are internally characterised by their local urban transport networks and 
at the same time are coupled through a long-distance inter-city transport network (like highways 
or flights). In our framework, this class of network can be induced by identifying the master 
graph $\bGamma$ with the macro-network, and the $M$ intra-layer subgraphs with the micro-networks. 
To have all micro-networks non-empty, the master graph must have all self-loops realised. This 
case is similar to the block-model mentioned above, but now the master graph itself can be 
chosen to have nontrivial structural properties, such as community structure, to resemble 
the specific properties of real-world networks of networks.

If the $t$-targeted degree sequences $\{\vec{k}_{s\to t}^*\}$ ($s,t=1,M$) are all enforced as constraints, 
then the relative entropy is given by \eqref{sMG} with $\gamma_{s,s}(\bGamma)=1$ for all $s$. 
However, in this class of models it is often more natural to assume that the internal degree sequence 
$\vec{k}_{s\to s}^*$ of each micro-network (layer) $s$ is enforced (in order to get realistic micro-network 
topologies), while between every pair $s,t$ ($s\ne t$) of micro-networks only the number of links 
$L^*_{s,t}$ is fixed (because the topology of the master graph is already chosen in order to replicate 
the empirical macro-network). This leads to the relaxed model in Theorem~\ref{resultcommrelax} 
with ${\cal D}=\left\lbrace (s,s)\colon\,s=1,M\right\rbrace$. The relative entropy is therefore
\be
s_\infty=\sum_{s=1}^M A_s\,\|f_{s\to s}\|_{\ell^1(g)}
\ee
and is still positive, even though the links among micro-networks do not contribute to it.


\section{Proof of the theorems}
\label{S3}


\subsection{Proof of Theorem \ref{resultunip}}
\label{S2.1}

\begin{proof}
The microcanonical number $\Omega_{\vec{k}^\star}$ is not known in general, but asymptotic 
results exist in the \emph{sparse regime} defined by the condition \eqref{eq:sparse1}. For this 
regime it was shown in~\cite{bender1977}, \cite{mckay1991} that
\begin{equation} 
\Omega_{\vec{k}^*} = \frac{\sqrt{2}\,(\frac{2L^*}{e})^{L^*}}{\prod_{i=1}^n k^*_i !}
\,e^{ - (\overline{k^{* 2}}/2\overline{k^*})^2 + \frac{1}{4} + o(n^{-1}\overline{k^*}^{\,3}) }, 
\label{peq:omegaCM}
\end{equation}
where $\overline{k^*} = n^{-1} \sum_{i=1}^n k^*_i$ (average degree), $L^*=n\overline{k^*}
/{2}$ (number of links), $\overline{k^{* 2}}= n^{-1} \sum_{i=1}^n k^{* 2}_i$ (average 
square degree). The canonical ensemble has Hamiltonian $H(\bG,\vec{\theta})
=\sum_{i=1}^n \theta_i k_i(\bG)$, where $\bG$ is a graph belonging to $\cG_{n}$, and 
$k_i(\bG) = \sum_{j\neq i} g_{i,j}(\bG)$ is the degree of the node $i$. The partition 
function equals 
\be
\begin{aligned}
Z(\theta) &= \sum_{\bG\in\cG_{n}} e^{-H(\bG,\vec{\theta})}
= \sum_{\bG\in\cG_{n}}\prod_{1\le i<j\le n} e^{-\theta_i g_{i,j}(\bG)}\\
&= \sum_{\bG\in\cG_{n}}\prod_{1\le i<j\le n} e^{-(\theta_i+\theta_j) g_{i,j}(\bG)}
=\prod_{1\le i<j\le n} (1+e^{-(\theta_i+\theta_j)}).
\end{aligned}
\ee
The canonical probability equals
\be
\Pcan(\bG \mid \vt) = \frac{\prod_{1\le i<j\le n} e^{-(\theta_i+\theta_j) g_{i,j}(\bG)}}{Z(\vt)}
=\prod_{1\le i<j\le n}\frac{ e^{-(\theta_i+\theta_j) g_{i,j}(\bG)}}{1+e^{-(\theta_i+\theta_j)}}.
\ee
Setting $p^*_{ij}\equiv e^{-\theta^*_i-\theta^*_j}/(1+e^{-\theta^*_i-\theta^*_j})$, and 
$\vec{\theta}^*$ such that 
\begin{equation}
\sum_{j\ne i} \frac{e^{-\theta^*_i-\theta^*_j}}{1+e^{-\theta^*_i-\theta^*_j}}=k^*_i
\quad\forall\,i
\label{peq:CM}
\end{equation}
we have
\begin{equation}
P_{\mathrm{can}}(\mathbf{G}) 
= \prod_{1\le i<j\le n} (p^*_{ij})^{g_{ij}} (1 - p^*_{ij})^{1 - g_{ij} }.
\label{peq:PCCM}
\end{equation}
It is ensured by \eqref{eq:sparse1} that $\lim_{n\to\infty} \frac{1}{n} \sum_{1\leq i<j\leq n}
{\widehat{p}_{ij}^{\,*\,2}} = 0$, a  condition under which we can show that \eqref{peq:PCCM} 
has the same asymptotic behaviour as 
\begin{equation}
\widehat{P}_{\mathrm{can}}(\mathbf{G}) 
= \prod_{1\le i<j\le n} (\widehat{p}^{\,*}_{ij})^{g_{ij}} (1 - \widehat{p}^{\,*}_{ij})^{1 - g_{ij} },
\end{equation}
with
\begin{equation}
\widehat{p}^{\,*}_{ij} = e^{-\theta^*_i-\theta^*_j}=\frac{k^*_i k^*_j}{2L^*}.
\label{peq:pijCL}
\end{equation}
Indeed,
\begin{equation}
\frac{1}{n}\log\left(\frac{\widehat{P}_{\mathrm{can}}(\mathbf{G}) }
{{P}_{\mathrm{can}}(\mathbf{G}) }\right)
= \frac{1}{n}\sum_{1\le i<j\le n}g_{i,j}\log(1-\widehat{p}_{ij}^{\,*})
-\frac{1}{n}\sum_{1\le i<j\le n}\log(1-\widehat{p}{_{ij}^{\,*}}^2) \to 0, \qquad n\to\infty,
\end{equation}
because 
\be
\sum_{1\leq i<j\leq n} g_{i,j} \log(1-\widehat{p}_{ij}^{\,*}) 
\leq (m^*)^2+O(\widehat{p}{_{ij}^{\,*}}^2)
\ee
and 
\be
0 \leq \frac{1}{n} \sum_{1\le i<j\le n}{\widehat{p}^{\,2}_{ij}}
=\frac{1}{2}\left[ \frac{\sum_{i=1}^n k_i^2}{\sqrt{n}\sum_{i=1}^n k_i}\right]^2
\leq \frac{1}{2}\frac{(m^*)^2}{{n}} \to 0, \qquad n\to\infty.
\ee
This implies $\sum_{1\le i<j\le n}\ln(1-p^*_{ij}) = - \sum_{1\leq i<j\leq n}k^*_i k^*_j/2L^*+o(n)$. Thus,
\begin{eqnarray}
\ln P_{\mathrm{can}}(\mathbf{G}^*) 
&=& \sum_{i = 1}^n k^*_i \ln k^*_i - L^* \ln (2L^*) - L^*+o(n). 
\label{peq:PCCL}
\end{eqnarray}
Combining \eqref{peq:omegaCM} and \eqref{peq:PCCL}, we obtain (recall \eqref{gdef})
\begin{equation}
\label{unipentropy}
S_n(\Pmic \mid \Pcan) = \sum_{i=1}^n g(k^*_i)+o(n), \qquad n\to\infty,
\end{equation}
where $g(k) = \log\left(\frac{k!}{k^k e^{-k}}\right)$, as defined in \eqref{gdef}. With the help of 
\eqref{empdef} this reads
\begin{equation}
\label{unipentropyemp}
n^{-1}\,S_n(\Pmic \mid \Pcan) = \sum_{k\in\N_0} f_n(k)g(k)+o(1) = \|f_n\|_{\ell^1(g)} +o(1),
\end{equation}
which together with \eqref{eq:conv1} yields the claim.
\end{proof}

\subsection{Proof of Theorem \ref{biplinks}}
\label{S2.22}

\begin{proof}
The microcanonical ensemble is easy: the number of graphs with a fixed 
fraction $\lambda \in (0,1)$ of links is 
\be
\Omega_{L^*} = {{n\choose 2} \choose L^*} = {K\choose \lambda K},\qquad K={n\choose 2}.
\ee
The canonical ensemble has the Hamiltonian $H(\bG,\theta)=\theta L(\bG)$, where 
$\bG$ is a graph belonging to $\cG_{n}$, and $L(\bG) =  
\sum_{1\le i<j\le n} g_{i,j}(\bG)$ is the number of links in $\bG$. The partition function equals 
\be
Z(\theta) = \sum_{\bG\in\cG_{n}} e^{-H(\bG,\theta)}
= \sum_{\bG\in\cG_{n}} \prod_{1\le i<j\le n} e^{-\theta g_{i,j}(\cG)}
=\prod_{1\le i<j\le n}(1+e^{-\theta}).
\ee
The canonical probability equals
\be
\Pcan(\bG \mid \theta) = \frac{e^{-\sum_{1\le i<j\le n} \theta g_{i,j}(\bG)}}{Z(\theta)}
=\prod_{1\le i<j\le n}\frac{e^{-\theta g_{i,j}(\bG)}}{1+e^{-\theta}}
=\prod_{1\le i<j\le n} p^{g_{i,j}(\bG)} (1-p)^{1-g_{i,j}(\bG)}
\ee
with $p=\frac{e^{-\theta}}{1+e^{-\theta}}$. We search for $\theta^*$ such that   
\be
L^*=\sum_{1\le i<j\le n} p^*, \qquad p^*=\frac{e^{-\theta^*}}{1+e^{-\theta^*}}.
\ee 
It follows that $p^*=\lambda$. Thus, 
\be
\begin{aligned}
\log {\Pmic(\bG^*)} &= -\log (K)! + \log (\lambda K)! + \log ((1-\lambda)K)!\\[0.2cm]
&=-K[\log K-1]+\lambda K[\log \lambda K -1]
+[(1-\lambda) K][\log((1-\lambda)K)-1]+o(n)\\
&= K\log(1-\lambda) + \lambda K\log\left( \frac{\lambda}{1-\lambda}\right)+o(n),\\
\log {\Pcan(\bG^*)} &= \lambda K\log(\lambda)+(1-\lambda)K\log(1-\lambda).
\end{aligned}
\ee
This in turn implies that
\be
\lim_{n\to \infty} \frac{S_n(\Pmic \mid \Pcan)}{n}
= 0.
\ee
\end{proof}

\subsection{Proof of Theorem \ref{result}}
\label{S2.4}

\begin{proof}
We start by describing the canonical ensemble. The Hamiltonian is
\be
\begin{aligned}
H(\bG|\vt,\vec{\phi}) &=\sum_{i\in\Lambda_1} k_i(\bG)\theta_i + \sum_{j\in\Lambda_2} k'_j(\bG)\phi_j\\
&= \sum_{i\in\Lambda_1} \sum_{j\in\Lambda_2}\theta_i g_{i,j}(\bG) 
+ \sum_{i\in\Lambda_1}\sum_{j\in\Lambda_2}\phi_j g_{i,j}(\bG)
=\sum_{i\in\Lambda_1}\sum_{j\in\Lambda_2}(\theta_i+\phi_j)g_{i,j}(\bG).
\end{aligned}
\ee
The partition function is
\be
Z(\vt,\vec{\phi}) = \sum_{\bG\in\cG_{{n_1},{n_2}}} e^{-\sum_{i\in\Lambda_1} \sum_{j\in\Lambda_2}
(\theta_i+\phi_j)g_{i,j}(\bG)}
= \prod_{i\in\Lambda_1} \prod_{j\in\Lambda_2} \left( 1+e^{-(\theta_i+\phi_j)} \right).
\ee
The canonical probability becomes
\be
\begin{aligned}
\Pcan(\bG \mid \vt,\vec{\phi}) 
&= \frac{e^{-\sum_{i\in\Lambda_1} \sum_{j\in\Lambda_2}(\theta_i+\phi_j) 
g_{i,j}(\bG)}}{Z(\vt,\vec{\phi})}\\
&= \prod_{i\in\Lambda_1} \prod_{j\in\Lambda_2} \frac{e^{-(\theta_i +\phi_j)
g_{i,j}(\bG)}}{1+e^{-(\theta_i+\phi_j)}}
=\prod_{i\in\Lambda_1}\prod_{j\in\Lambda_2}p_{i,j}^{g_{i,j}(\bG)}(1-p_{i,j})^{1-g_{i,j}(\bG)},
\end{aligned}
\ee
where $p_{i,j}=\frac{e^{-(\theta_i +\phi_j)}}{1+e^{-(\theta_i+\phi_j)}}$. We search for 
$(\vt^*,\vec{\phi}^*)$ that solves the system of equations
\be
\begin{cases}
\sum_{j\in\Lambda_2} p_{i,j}^*=k_i^*,\\
\sum_{i\in\Lambda_1} p_{i,j}^*=k'^*_j,
\end{cases}
\ee
where $p_{i,j}^*=\frac{e^{-(\theta_i^* +\phi_j^*)}}{1+e^{-(\theta_i^*+\phi_j^*)}}$. If $\mathbf{G^*}$ 
is any graph in $\cG_{{n_1},{n_2}}$ such that $k_i(\mathbf{G^*})=k_i^*$ and $k'_j(\mathbf{G^*})=k'^*_j$, 
then 
\begin{equation}
\label{canonicalbip}
\Pcan(\bG)=\prod_{i\in\Lambda_1}\prod_{j\in\Lambda_2} {p_{i,j}^*}^{g_{i,j}(\bG)}(1-p^*_{i,j})^{1-g_{i,j}(\bG)}.
\end{equation}

Under the sparseness condition \eqref{maxdegree}, we can replace $p_{i,j}^*$ with the following quantity. Define $\widehat{p}_{i,j}^{\,*}
= e^{-(\theta_i^* +\phi_j^*)}$ and consider the system of equations
\begin{equation}
\begin{cases}
\sum_{j\in\Lambda_2} \widehat{p}_{i,j}^{\,*}=k^*_i,\\
\sum_{i\in\Lambda_1} \widehat{p}_{i,j}^{\,*}=k'^*_j.
\end{cases}
\end{equation}
This has solution 
\be
\widehat{p}_{i,j}^{\,*}=\frac{k^*_i k'^*_j}{L^*}, \qquad 
L^*=\sum_{i\in \Lambda_1} k^*_i=\sum_{j\in \Lambda_2} k'^*_j.
\ee
We define 
\begin{equation}
\widehat{P}_{\mathrm{can}}(\mathbf{G}) 
= \prod_{i\in\Lambda_1}\prod_{j\in\Lambda_2} (\widehat{p}^{\,*}_{i,j})^{g_{ij}(\bG)} 
(1 - \widehat{p}^{\,*}_{i,j})^{1 - g_{ij}(\bG)},
\end{equation}
and note that 
\be
\frac{1}{n_1+n_2}\log\left(\frac{\widehat{P}_{\mathrm{can}}
(\mathbf{G}) }{{P}_{\mathrm{can}}(\mathbf{G}) } \right) \to 0, \qquad {n_1},{n_2}\to\infty.
\ee
The crucial point is to prove that $\frac{1}{{n_1}+{n_2}}\sum_{i\in\Lambda_1}\sum_{j\in\Lambda_2} 
\widehat{p}^{\,*\,2}_{i,j}\to 0$. This allows us to write  
\begin{equation}
\label{appr}
\sum_{i\in\Lambda_1}\sum_{j\in\Lambda_2}\log(1-p^*_{i,j}) = 
-\sum_{i\in\Lambda_1}\sum_{j\in\Lambda_2}\frac{k_i^*k'^*_j}{L^*}+o({n_1}+{n_2}), 
\qquad {n_1},{n_2}\to \infty.
\end{equation}
Indeed, 
\be
0\leq \frac{1}{{n_1}+{n_2}}\sum_{i\in\Lambda_1}\sum_{j\in\Lambda_2} \widehat{p}^{\,*\,2}_{i,j}
=\frac{1}{{n_1}+{n_2}}\frac{\sum_{i\in\Lambda_1} {k^*_i}^2 
\sum_{j\in\Lambda_2}{k'^*_j}^2}{\sum_{i\in\Lambda_1} k^*_i\sum_{j\in\Lambda_2} k'^*_j}
\leq \frac{m^*m'^*}{\sqrt{{n_1}{n_2}}}\frac{\sqrt{{n_1}{n_2}}}{{n_1}+{n_2}}\to 0,
\ee
because ${m^*m'^*}=o({L^*}^{2/3})$ implies $m^*m'^*=o(\sqrt{{n_1}{n_2}})$. 

Combining \eqref{canonicalbip} and \eqref{appr}, we have
\begin{eqnarray}
\label{logcan}
\log \Pcan(\mathbf{G^*}) 
= \sum_{i\in\Lambda_1}\sum_{j\in\Lambda_2} g_{i,j}(\bG^*)\log\left(\frac{k_i^*k'^*_j}{L^*}\right)
- \sum_{i\in\Lambda_1}\sum_{j\in\Lambda_2}\frac{k_i^*k'^*_j}{L^*}+o({n_1}+{n_2})\nonumber\\
= \sum_{i\in\Lambda_1} k_i^*\log\left(k_i^*\right) + \sum_{j\in\Lambda_2} k'^*_j
\log\left(k'^*_j\right)-L^*\log L^* -L^*+o({n_1}+{n_2}),
\end{eqnarray}
which concludes our computation for the canonical ensemble.

Microcanonical probabilities come from the results in \cite{GMW06}, where it is shown 
that, as $n\to\infty$, the number of bi-partite graphs with degree sequences 
$\vec{k^*},\vec{{k'}^*}$ on the two layers is given by 
\be
\Omega_{\vec{k^*},\vec{{k'}^*}} = \frac{L^*!}{\prod_{i\in\Lambda_1} k_i^*! 
\prod_{j\in\Lambda_2} k'^*_j!}\,e^{o({n_1}+{n_2})}.
\ee
Hence
\be
\label{logMicroc}
\log \Pmic(\mathbf{G^*}) = - \log \Omega_{\vec{s^*},\vec{t^*}}
=\sum_{i\in\Lambda_1} k^*_i!+\sum_{j\in\Lambda_2} {k'}^*_j!-\log (L^*!)+o({n_1}+{n_2}).
\ee
From \eqref{logcan} and \eqref{logMicroc} we get
\be
\label{relativebipartite}
\begin{aligned}
S_{{n_1}+{n_2}}(\Pcan \mid \Pmic) &= \log \Pmic(\mathbf{G^*}) - \log \Pcan(\mathbf{G^*})\\
& = \sum_{i\in\Lambda_1} \log \left(\frac{k_i^* !}{{k_i^*}^{k_i^*}}\right)
+ \sum_{j\in\Lambda_2} \log \left(\frac{k'^*_j !}{{k'^*_j}^{k'^*_j}}\right)\\
&\qquad +L^*\log L^* +L^*-\log (L^*!)+o({n_1}+{n_2})\\
&=\sum_{i\in\Lambda_1} g(k_i^*)+\sum_{j\in\Lambda_2} g(k'^*_j)+o({n_1}+{n_2}),
\end{aligned}
\ee
where in the last line we use $L^*=\sum_{i\in \Lambda_1} k^*_i=\sum_{j\in \Lambda_2} {k'}^*_j$ 
and Stirling's approximation for $\log(L^* !)$. Since
\be
\begin{aligned}
{n_1}^{-1} \sum_{i\in\Lambda_1} g(k_i^*) &= \sum_{k\in\N_0}^{n_2} f_{1\to 2}^{({n_1})}(k) g(k)
= \|f_{1\to 2}^{({n_1})}\|_{\ell^1(g)},\\
{n_2}^{-1} \sum_{j\in\Lambda_2} g(k'^*_j) &= \sum_{k\in\N_0}^{n_1} f_{2\to 1}^{({n_2})}(k) g(k)
= \|f_{2\to 1}^{({n_2})}\|_{\ell^1(g)},
\end{aligned}
\ee
we get, with the help of \eqref{convergence},
\be
\lim_{n\to\infty} \frac{S_{{n_1}+{n_2}}(\Pcan \mid \Pmic)}{{n_1}+{n_2}}
= A_1\,\|f_{1\to 2}\|_{\ell^1(g)} + A_2\,\|f_{2\to 1}\|_{\ell^1(g)},
\ee
which proves the claim.
\end{proof}


\subsection{Proof of Theorem \ref{configurationmodel}}
\label{S2.3}

\begin{proof}
The number of bi-partite graphs with constraint ${\vec{k}^*}$ on the top layer is
\begin{equation}
\Omega_{\vec{k}^*} = \prod_{i\in\Lambda_1} {{n_2} \choose k_i^*}. 
\end{equation}
In order to calculate the canonical probability, we calculate the partition function:
\be
Z(\vt) = \sum_{\bG\in\cG_{{n_1},{n_2}}} e^{-\sum_{i\in \Lambda_1} 
\theta_i \sum_{j\in\Lambda_2} g_{i,j}(\bG)}
= \sum_{\bG\in\cG_{{n_1},{n_2}}} \prod_{i\in\Lambda_1} \prod_{j\in \Lambda_2} e^{-\theta_i g_{i,j}(\bG)}
= \prod_{i\in\Lambda_1} \prod_{j\in\Lambda_2} [1+e^{-\theta_i}].
\ee
The canonical probability becomes
\be
\Pcan(\bG|\vt) = \frac{e^{-\sum_{i\in\Lambda_1} \theta_i \sum_{j\in\Lambda_2} g_{i,j}(\bG)}}{Z(\vt)}
= \prod_{i\in\Lambda_1} \prod_{j\in\Lambda_2} \frac{e^{-\theta_i g_{i,j}(\bG)}}{1+e^{-\theta_i}}
= \prod_{i\in\Lambda_1} \prod_{j\in\Lambda_2} p_i^{g_{i,j}(\bG)} (1-p_i)^{1-g_{i,j}(\bG)}
\ee
with $p_i=\frac{e^{-\theta_i}}{1+e^{-\theta_i}}$. We search for $\theta^*_i$ such that   
\be
k_i^* = \sum_{j\in \Lambda_2} p^*_{i} = {n_2}p^*_i, \qquad p^*_i=\frac{e^{-\theta^*_i}}{1+e^{-\theta^*_i}}.
\ee
It follows that $p_i=\frac{k_i^*}{n_2}$ (recall \eqref{eq:PCRG}). According to \eqref{eq:KL2} 
we have 
\be
\begin{aligned}
S_{{n_1}+{n_2}}(\Pmic \mid \Pcan) &= \ln \frac{\Pmic(\bG^*)}{\Pcan(\bG^*)}\\
&= -\sum_{i\in\Lambda_1}\log{{n_2}\choose k_i^*} 
- \sum_{i\in\Lambda_1} k_i^*\log\left(\dfrac{k_i^*}{n_2}\right)
-\sum_{i\in\Lambda_1} ({n_2}-k_i^*)\log\left(1-\dfrac{k_i^*}{n_2}\right)\\
&= {n_1}{n_2}\log {n_2} -\sum_{i\in\Lambda_1}  
\log\left[ {{n_2} \choose k_i^*}{k_i^*}^{k_i^*}{({n_2}-k_i^*)}^{({n_2}-k_i^*)}\right].
\end{aligned}
\ee
Abbreviate $U_{a}(x)\equiv\log\left[ {a \choose x}{x}^{x}{(a-x)}^{a-x}\right]$ and write
\be
S_{{n_1}+{n_2}}(\Pmic \mid \Pcan) 
= {n_1}{n_2}\log {n_2} - \sum_{i\in\Lambda_1}  U_{n_2}(k_i^*)
= {n_1}{n_2}\log {n_2} -{n_1} \sum_{k=0}^{n_2} f_{n_1}(k)U_{n_2}(k).
\ee
For the relative entropy per node this gives
\be
\begin{aligned}
s_{{n_1}+{n_2}} &= \frac{n_1}{{n_1}+{n_2}} \sum_{k=0}^{n_2} f_{n_1}(k)\,{n_2}\log {n_2}
- \frac{n_1}{{n_1}+{n_2}}\sum_{k=0}^{n_2}f_{n_1}(k) U_{n_2}(k)\\
&= -\frac{n_1}{n_1+{n_2}} \sum_{k=0}^{n_2} f_{n_1}(k) 
\log \mathrm{Bin}\big({n_2},\tfrac{k}{n_2}\big)(k)=\frac{n_1}{{n_1}+{n_2}}\|f_{n_1}\|_{\ell^1(g_{n_2})}.
\end{aligned}
\ee


\paragraph{Case (1).}
Recall \eqref{gdef}. Note that $x \mapsto z(x)= e^{g(x)}$ is non-decreasing:
\be
\frac{z(x-1)}{z(x)} = \left({\frac{x}{x-1}}\right)^{x-1} \frac{1}{e} \leq 1.
\ee
It therefore follows that
\be
\label{uest1}
\begin{aligned}
&\|f_{n_1}\|_{\ell^1(g_{n_2})}=-\sum_{k=0}^{n_2} f_{n_1}(k) 
\log \mathrm{Bin}\big({n_2},\tfrac{k}{n_2}\big)(k)
= \sum_{k=0}^{n_2} f_{n_1}(k) \log\left(\frac{z(k)z({n_2}-k)}{z({n_2})}\right)\\
&= \sum_{k\in\N_0} f_{n_1}(k) \log\left(\frac{z(k)z({n_2}-k)}{z({n_2})}\right) \mathbb{I}_{k \leq {n_2}}
\leq \sum_{k\in\N_0} \mathbb{I}_{0\leq k \leq {n_2}}\,f_{n_1}(k)\log z(k)\\
&\leq \sum_{k\in\N_0} f_{n_1}(k) \log z(k) = \|f_{n_1}\|_{\ell^1(g)} < \infty.
\end{aligned}
\ee
By \eqref{twofixedcon} and dominated convergence, we may exchange limit and 
sum to obtain
\be
\lim_{{n}\to\infty} s_{{n_1},{n_2}}
= \lim_{{n_2}\to\infty} \frac{n_1}{{n_1}+{n_2}} \sum_{k\in\N_0} f_{n_1}(k)
\lim_{{n_2}\to\infty} \log\left(\frac{z(k)z({n_2}-k)}{z({n_2})}\right) \mathbb{I}_{0\leq k \leq {n_2}} = 0,
\ee
where we use that $\lim_{{n}\to\infty}\frac{n_1}{{n_1}+{n_2}}=0$ and $\lim_{{n}\to\infty} 
\frac{z({n_2}-k)}{z({n_2})}=1$ for all $k\in\N_0$.  


\paragraph{Case (2).}
Using \eqref{uest1} and \eqref{twofixedcon}, we get
\be
0\le s_{{n_1},{n_2}}
= \frac{n_1}{{n_1}+{n_2}} \|f_{n_1}\|_{\ell^1(g)} \to^{n\to\infty} \frac{1}{1+c}\|f\|_{\ell^1(g)}=0.
\ee


\paragraph{Case (3).}
Estimate
\be
0\le |\|f_{n_1}\|_{\ell^1(g_{n_2})}-\|f_{n_1}\|_{\ell^1(g_{n_2})}|
\le \|f_{n_1}-f\|_{\ell^1(g_{n_2})} \le \|f_{n_1}-f\|_{\ell^1(g)}\to^{{n}\to \infty} 0.
\ee


\paragraph{Case (4).}
\be
0\le |\|f_{n_1}\|_{\ell^1(g_{n_2})}-\|f\|_{\ell^1(g)}|\le
\sum_{k\in \mathbb{N}_0}|f_{n_1}(k)-f(k)||g_{n_2}(k)\mathbb{I}_{k \leq {n_2}}-g(k)|
\le 2\|f_{n_1}-f\|_{\ell^1(g)}.
\ee
Since $\frac{n_1}{{n_1}+{n_2}} = \frac{1}{1+\frac{n_2}{n_1}}
\to\frac{1}{1+c}$, the claim follows.
\end{proof}


\subsection{Proof of Theorem \ref{resultbip}}
\label{S2.2}

\begin{proof}
The microcanonical ensemble is easy: the number of bi-partite graphs with a fixed 
fraction $\lambda \in (0,1)$ of links is 
\be
\Omega_{L^*} = {{n_1}{n_2} \choose L^*} = {{n_1}{n_2} \choose \lambda {n_1}{n_2}}.
\ee
The canonical ensemble has the Hamiltonian $H(\bG,\theta)=\theta L(\bG)$, where 
$\bG$ is a bi-partite graph belonging to $\cG_{{n_1},{n_2}}$, and $L(\bG) = \sum_{i\in\Lambda_1}
\sum_{j\in\Lambda_2} g_{i,j}(\bG)$ is the number of links in $\bG$. The partition function equals 
\be
Z(\theta) = \sum_{\bG\in\cG_{{n_1},{n_2}}} e^{-H(\bG,\theta)}
= \sum_{\bG\in\cG_{{n_1},{n_2}}} \prod_{i\in\Lambda_1} \prod_{j\in\Lambda_2} e^{-\theta g_{i,j}(\cG)}
=\prod_{i\in\Lambda_1} \prod_{j\in\Lambda_2} (1+e^{-\theta}).
\ee
The canonical probability equals
\be
\Pcan(\bG \mid \vt) = \frac{e^{-\sum_{i\in\Lambda_1}  \sum_{j\in\Lambda_2} \theta g_{i,j}(\bG)}}{Z(\vt)}
=\prod_{i\in\Lambda_1} \prod_{j\in\Lambda_2} \frac{e^{-\theta g_{i,j}(\bG)}}{1+e^{-\theta}}
=\prod_{i\in\Lambda_1} \prod_{j\in\Lambda_2} p^{g_{i,j}(\bG)} (1-p)^{1-g_{i,j}(\bG)}
\ee
with $p=\frac{e^{-\theta}}{1+e^{-\theta}}$. We search for $\theta^*$ such that   
\be
L^*=\sum_{i\in\Lambda_1}\sum_{j\in\Lambda_2} p^*, \qquad p^*=\frac{e^{-\theta^*}}{1+e^{-\theta^*}}.
\ee 
It follows that $p^*=\lambda$. Thus, 
\be
\begin{aligned}
\log {\Pmic(\bG^*)} &= -\log ({n_1}{n_2})! + \log (\lambda {n_1}{n_2})! + \log ((1-\lambda){n_1}{n_2})!\\
&= -{n_1}{n_2}[\log {n_1}{n_2}-1]+\lambda {n_1}{n_2}[\log \lambda {n_1}{n_2} -1]\\
&\qquad +[1-\lambda {n_1}{n_2}][\log((1-\lambda){n_1}{n_2})-1]+o({n_1}+{n_2})\\
&= {n_1}{n_2}\log(1-\lambda) + \lambda {n_1}{n_2}
\log\left( \frac{\lambda}{1-\lambda}\right)+o({n_1}+{n_2}),\\
\log {\Pcan(\bG^*)} &= {n_1}{n_2}\log(1-\lambda) 
+ \lambda {n_1}{n_2}\log\left( \frac{\lambda}{1-\lambda}\right).
\end{aligned}
\ee
This in turn implies that
\be
\lim_{{n_1},{n_2}\to \infty} \frac{S_{{n_1}+{n_2}}(\Pmic \mid \Pcan)}{{n_1}+{n_2}}= 0.
\ee
\end{proof}

\subsection{Proof of Theorem \ref{resultcomm}}
\label{S2.5}

\begin{proof}
The proof is based on the previous theorems. We start by looking at the Hamiltonian of the system. 
For each admitted pair of layers ($\gamma_{s,t}(\bGamma)=1$) we define Lagrange multipliers 
$\vt_{s\to t} = (\theta_{1}^{(t)},\dots,\theta_{{n_s}}^{(t)})$. The Hamiltonian equals 
\be
\begin{aligned}
&H\big(\bG \mid \vt_{s\to t};\ s,t=1,\dots,M,\ \gamma_{s,t}(\bGamma)=1\big)\\ 
&= \sum_{\substack{1\le s<t\le M\\ \gamma_{s,t}(\bGamma)=1}}
\sum_{\substack{i\in \Lambda_s\\j\in  \Lambda_t}}(\theta_i^t+\theta_j^s)g_{i,j}(\bG)
+\sum_{\substack{s=1\\ \gamma_{s,s}(\bGamma)=1}}^M 
\sum_{\substack{i,j\in \Lambda_s\\i<j}} (\theta_i^s+\theta_j^s) g_{i,j}(\bG)\\
&= \sum_{\substack{1\le s<t\le M\\ \gamma_{s,t}(\bGamma)=1}}
\sum_{\substack{i\in \Lambda_s\\j\in  \Lambda_t}} H_{s,t} (\bG^{(st)} \mid \vt_{s\to t},\vt_{t\to s})
+\sum_{\substack{s=1\\ \gamma_{s,s}(\bGamma)=1}}^M 
\sum_{\substack{i,j\in \Lambda_s\\i<j}} H_{s,s}(\bG^{(ss)} \mid\vt_{s\to s}),
\end{aligned}
\ee
where 
\be
\begin{aligned}
H_{s,t}(\bG^{(st)} \mid \vt_{s\to t},\vt_{t\to s})
&= \sum_{\substack{i\in \Lambda_s\\j\in  \Lambda_t}} (\theta_i^t+\theta_j^s) g_{i,j}(\bG^{(st)}),\\
H_{s,s}(\bG^{(ss)} \mid \vt_{s\to s})
&= \sum_{\substack{i,j\in \Lambda_s\\i<j}}(\theta_i^s+\theta_j^s)g_{i,j}(\bG^{(ss)}),
\end{aligned}
\ee
and $\bG^{(st)}$ ($\bG^{(ss)}$) is the bi-partite (uni-partite) graph between layers $s$ and $t$ 
(inside layer $s$) obtained from the multi-partite graph $\bG$. The $n_s\times n_t$ matrix 
representing the bi-partite graph has, for each $i\in \Lambda_s$ and $j\in \Lambda_t$, elements 
$g_{i,j}(\bG^{(st)})=g_{i,j}(\bG)$. Note that $H_{s,t}(\bG^{(st)}\mid\vt_{s\to t},\vt_{t\to s})$ is the 
Hamiltonian of the bi-partite graph $\bG^{(st)}$ between layers $s$ and $t$ with constraints 
${\vec{k}}^{\,*}_{s\to t}$, and $H_{s,s}(\bG^{(ss)} \mid \vt_{s\to s})$ is the Hamiltonian of the 
uni-partite graph $\bG^{(ss)}$ of the layer $s$ with constraints ${\vec{k}}^{\,*}_{s\to s}$.

The partition function of the canonical ensemble equals
\be
\begin{aligned}
&Z\big(\vt_{s\to t};\ s,t=1,\dots,M,\ \gamma_{s,t}(\bGamma)=1\big)
= \sum_{\bG\in\cG_{n_1,\dots ,{n_M}}(\bGamma)}
e^{-H(\bG\,\mid\,\vt_{s\to t};\ s,t=1,2,\dots,M\colon\, \gamma_{s,t}(\bGamma)=1)}\\
&= \prod_{\substack{1\le s<t\le M\\ \gamma_{s,t}(\bGamma)=1}}
\sum_{\bG^{(st)}\in\cG_{n_s,{n_t}}}e^{-H_{s,t}(\bG^{(st)} \,\mid\, \vt_{s\to t},\vt_{t\to s})}
\prod_{\substack{s=1\\ \gamma_{s,s}(\bGamma)=1}}^M\sum_{\bG^{(ss)}\in\cG_{n_s,{n_s}}}
e^{-H_{s,s}(\bG^{(ss)}\,\mid\, \vt_{s\to s})}\\
&= \prod_{\substack{1\le s<t\le M\\ \gamma_{s,t}(\bGamma)=1}}
Z^{(st)}(\vt_{s\to t},\vt_{t\to s})\prod_{\substack{s=1\\ \gamma_{s,s}(\bGamma)=1}}^M 
Z^{(ss)}(\vt_{s\to s}),
\end{aligned}
\ee
where $Z^{(st)}(\vt_{s\to t},\vt_{t\to s})$ is the partition function of the set of bi-partite 
graphs $\cG_{n_s,n_t}$ with constraints ${\vec{k}}^{\,*}_{s\to t}$ on the top layer and 
${\vec{k}}^{\,*}_{t\to s}$ on the bottom layer, and $Z^{(ss)}(\vt_{s\to s})$ is the partition 
function of the set of graph $\cG_{n_s}$ with constraint ${\vec{k}}^{\,*}_{s\to s}$. The 
canonical ensemble is
\begin{eqnarray}
\label{cancomm}
\Pcan(\bG)=\prod_{\substack{1\le s<t\le M\\ \gamma_{s,t}(\bGamma)=1}}
P^{(st)}_{\mathrm{can}}(\bG^{(st)})
\prod_{\substack{s=1\\ \gamma_{s,s}(\bGamma)=1}}^M P^{(ss)}_{\mathrm{can}}(\bG^{(ss)}),
\end{eqnarray}
where $P^{(st)}_{\mathrm{can}}(\bG^{(st)})$ is the canonical probability of the bi-partite 
graph $\bG^{(st)}$ with constraints ${\vec{k}}^{\,*}_{s\to t}$ on the top layer and ${\vec{k}}^{\,*}_{t\to s}$ 
on the bottom layer, and $P^{(ss)}_{\mathrm{can}}(\bG^{(ss)})$ is the canonical probability 
of the uni-partite graph $ \bG^{(ss)} $ with constraint ${\vec{k}}^{\,*}_{s\to s}$.

We can split the microcanonical probability as products of microcanonical probabilities 
for simpler cases. The number of graphs with constraints $\vC^*$ is
\be
\begin{aligned}
&\Omega_{{\vec{k}}^{\,*}_{s\to t};\,s,t\in\left\lbrace1,\dots,M\right\rbrace,\ \gamma_{s,t}(\bGamma)=1} 
= \left|\left\{ \bG\in\cG_{n_1,\dots,n_{M}}(\bGamma)\colon\hspace{-.2cm} 
\sum_{j\in  \Lambda_t} g_{i,j}(\bG) = k_i^{*\,t}\,\,
\forall\,i\in \Lambda_s\,\forall\,\gamma_{s,t}=1\right\}\right|\\
&=\prod_{\substack{1\le s<t\le M\\ \gamma_{s,t}(\bGamma)=1}} 
\left|\left\{ \bG^{(st)}\in\cG_{n_s,n_{k}}\colon \hspace{-.2cm}\sum_{j\in  \Lambda_t} 
g_{i,j}(\bG^{(st)})=k_i^{*\,t}\,\,\forall\,i\in \Lambda_s,
\sum_{i\in \Lambda_s} g_{i,j}(\bG^{(st)})=k_j^{*\,s}\,\forall\,j\in \Lambda_t\right\}\right|\\
&\qquad \prod_{\substack{s=1\\ \gamma_{s,s}(\bGamma)=1}}^M 
\left|\left\{ \bG^{(ss)}\in\cG_{n_s}\colon\hspace{-.2cm} 
\sum_{j\in \Lambda_s}g_{i,j}(\bG^{(ss)})=s_i^{*\,h}\,\,\forall\,i\in \Lambda_s \right\}\right|\\
&=\prod_{\substack{1\le s<t\le M\\ \gamma_{s,t}(\bGamma)=1}} 
\Omega_{{\vec{k}}^{\,*}_{s\to t},{\vec{k}}^{\,*}_{t\to s}}
\prod_{\substack{s=1\\ \gamma_{s,s}(\bGamma)=1}}^M \Omega_{{\vec{k}}^{\,*}_{s\to s}}.
\end{aligned}
\ee
This means the microcanonical probability can be factorised as
\begin{eqnarray}
\label{miccomm}
\Pmic(\bG)=\prod_{\substack{1\le s<t\le M\\ \gamma_{s,t}(\bGamma)=1}}P^{(st)}_{\mathrm{mic}}(\bG^{(st)})
\prod_{\substack{s=1\\ \gamma_{s,s}(\bGamma)=1}}^M P^{(ss)}_{\mathrm{mic}}(\bG^{(ss)}),
\end{eqnarray}
where $P^{(st)}_{\mathrm{mic}}(\bG^{(st)})$ is the microcanonical probability of the bi-partite 
graph $\bG^{(st)}$ with constraints ${\vec{k}}^{\,*}_{s\to t}$ on the top layer and ${\vec{k}}^{\,*}_{t\to s}$ 
on the bottom layer, and  $P^{(ss)}_{\mathrm{mic}}(\bG^{(ss)})$ is the microcanonical probability 
of the uni-partite graph $\bG^{(ss)}$ with constraint ${\vec{k}}^{\,*}_{s\to s}$.

Equations \eqref{cancomm} and \eqref{miccomm} imply that the relative entropy equals the sum
\begin{eqnarray}
S_n(\Pmic \mid \Pcan) = \sum_{\substack{1\le s<t\le M\\ \gamma_{s,t}(\bGamma)=1}} 
S_n(P^{(st)}_{\mathrm{mic}} \mid P^{(st)}_{\mathrm{can}})
+\sum_{\substack{s=1\\ \gamma_{s,s}(\bGamma)=1}}^M 
S_n(P^{(ss)}_{\mathrm{mic}} \mid P^{(ss)}_{\mathrm{can}}).
\end{eqnarray}
We can now apply Theorems \ref{resultunip} and \ref{result} to get the asymptotic relative entropy 
per nodes as
\be
\begin{aligned}
&\lim_{n_1,\dots,n_{M}\to\infty}\frac{S_n(\Pmic \mid \Pcan)}{n}\\
&= \sum_{\substack{1\le s<t\le M\\ \gamma_{s,t}(\bGamma)=1}}\lim_{n_1,\dots,n_{M}\to\infty}
\frac{S_n(P^{(st)}_{\mathrm{mic}} \mid P^{(st)}_{\mathrm{can}})}{n}
+ \sum_{\substack{s=1\\ \gamma_{s,s}(\bGamma)=1}}^M\lim_{n_1,\dots,n_{M}\to\infty} 
\frac{S_n(P^{(ss)}_{\mathrm{mic}} \mid P^{(ss)}_{\mathrm{can}})}{n}\\
&= \sum_{\substack{1\le s<t\le M\\ \gamma_{s,t}(\bGamma)=1}}
\left\lbrace A_s\,\|f_{s\to t}\|_{\ell^1(g)} +A_t\,\|f_{t\to s}\|_{\ell^1(g)} \right\rbrace 
+\sum_{\substack{s=1\\ \gamma_{s,s}(\bGamma)=1}}^M \left\lbrace A_s\,
\|f_{s\to s}\|_{\ell^1(g)} \right\rbrace\\
&=\sum_{\substack{s,t=1\\ \gamma_{s,t}(\bGamma)}}^M A_s\,\| f_{s\to t}\|_{\ell^1(g)}.
\end{aligned}
\ee
\end{proof}


\subsection{Proof of Theorem \ref{resultcommrelax}}
\label{S2.6}

\begin{proof}
We start by studying the Hamiltonian. For each pair $(s,t)$ of layers in $\mathcal D$, we 
define Lagrange multipliers $\vt_{s\to t}=(\theta_{1}^t,\dots,\theta_{{n_s}}^t)$. 
For each pair $(s,t)$ of layers in $\mathcal L$, we define a Lagrange multiplier 
${\theta}_{s,t}$. The Hamiltonian is 
\be
\begin{aligned}
&H\big(\bG \mid \vt_{s\to t},\theta_{l,m};\, (s,t)\in \mathcal D, (l,m)\in \mathcal L\big)\\
&=H_{\mathcal D}(\bG \mid \vt_{s\to t};\, (s,t)\in \mathcal D)
+ H_{\mathcal L}(\bG \mid \theta_{l,m};\, (l,m)\in \mathcal L)
\end{aligned}
\ee
with 
\be
\begin{aligned}
H_{\mathcal D}(\bG \mid \vt_{s\to t};\, (s,t)\in \mathcal D)
&= \sum_{\substack{1\le s<t\le M\\(s,t)\in \mathcal D}}
\sum_{\substack{i\in \Lambda_s\\j\in  \Lambda_t}}
(\theta_i^t+\theta_j^s)g_{i,j}(\bG)
+\sum_{{\substack{s=1\\(s,s)\in \mathcal D}}}^M \sum_{\substack{i,j\in \Lambda_s\\i<j}}
(\theta_i^s+\theta_j^s)g_{i,j}(\bG),\\
H_{\mathcal L}(\bG \mid {\theta}_{s,t};\,(s,t)\in \mathcal L)
&=\sum_{\substack{1\le s<t\le M\\(s,t)\in \mathcal L}}
\sum_{\substack{i\in \Lambda_s\\j\in  \Lambda_t}}
(\theta_{s,t})g_{i,j}(\bG)
+\sum_{{\substack{s=1\\(s,s)\in \mathcal L}}}^M
\sum_{\substack{i,j\in \Lambda_s\\i<j}}
(\theta_{s,s})g_{i,j}(\bG).
\end{aligned}
\ee
Consequently, the canonical ensemble is
\begin{eqnarray}
\Pcan(\bG)=P^{\mathcal D}_{\mathrm{can}}(\bG)P^{\mathcal L}_{\mathrm{can}}(\bG)
\end{eqnarray}
with
\be
\begin{aligned}
P^{\mathcal D}_{\mathrm{can}}(\bG) 
&= \prod_{\substack{1\le s<t\le M\\(s,t)\in \mathcal D}}{P^{(st)}_{\mathrm{can}}}^{\mathcal D}(\bG^{(st)})
\prod_{\substack{s=1\\(s,s)\in \mathcal D}}^M {P^{(ss)}_{\mathrm{can}}}^{\mathcal D}(\bG^{(ss)}),\\
P^{\mathcal L}_{\mathrm{can}}(\bG)
&= \prod_{\substack{1\le s<t\le M\\(s,t)\in \mathcal L}}{P^{(st)}_{\mathrm{can}}}^{\mathcal L}(\bG^{(st)})
\prod_{\substack{s=1\\(s,s)\in \mathcal L}}^M {P^{(ss)}_{\mathrm{can}}}^{\mathcal L}(\bG^{(ss)}).
\end{aligned}
\ee
Here,
\begin{itemize} 
\item
$\bG^{(st)}$ ($\bG^{(ss)}$ ) is the bi-partite (uni-partite) graph between layers $s$ and $t$ (and itself) 
obtained from the multi-partite graph $\bG$. The $n_s\times n_t$ ($n_s\times n_s$) matrix representing 
this bi-partite (uni-partite) graph has, for each $i\in \Lambda_s$ and $j\in \Lambda_t$ (for each $i,j\in 
\Lambda_s$), elements $g_{i,j}(\bG^{(st)})=g_{i,j}(\bG)$ ($g_{i,j}(\bG^{(ss)})=g_{i,j}(\bG)$).
\item
${P^{(st)}_{\mathrm{can}}}^{\mathcal D}(\bG^{(st)})$ (${P^{(ss)}_{\mathrm{can}}}^{\mathcal D}(\bG^{(ss)})$) 
is the canonical probability of the bi-partite (uni-partite) graph $\bG^{(st)}$ ($\bG^{(ss)}$) with constraints 
${\vec{k}}^{\,*}_{s\to t}$ on the top layer and ${\vec{k}}^{\,*}_{t\to s}$ on the bottom layer (with constraint 
${\vec{k}}^{\,*}_{s\to s}$). 
\item
${P^{(st)}_{\mathrm{can}}}^{\mathcal L}(\bG^{(st)})$ (${P^{(ss)}_{\mathrm{can}}}^{\mathcal L}(\bG^{(ss)})$) 
is the canonical probability of the bi-partite (uni-partite) graph $\bG^{(st)}$ ($\bG^{(ss)}$) with constraint 
$L_{s,t}^*$ ($L^*_{s,s}$). 
\end{itemize}

We can split the microcanonical probability as products of microcanonical probabilities 
of simpler cases. The number of graphs with such a type of constraints is 
\begin{eqnarray}
\Omega_{{\vec{k}}^{\,*}_{s\to t}, L^*_{l,m};\,(s,t)\in \mathcal D, (l,m)\in \mathcal L}
=\Omega_{{\vec{k}}^{\,*}_{s\to t};(s,t)\in \mathcal D}\,\Omega_{L_{l,m};(l,m)\in \mathcal L}.
\end{eqnarray}
This means that the microcanonical probability can be factorised as
\begin{eqnarray}
\Pmic(\bG)=P^{\mathcal D}_{\mathrm{mic}}(\bG)P^{\mathcal L}_{\mathrm{mic}}(\bG)
\end{eqnarray}
with
\be
\begin{aligned}
&P^{\mathcal D}_{\mathrm{mic}}(\bG)
=\prod_{\substack{1\le s<t\le M\\(s,t)\in \mathcal D}}{P^{(st)}_{\mathrm{mic}}}^{\mathcal D}(\bG^{(st)})
\prod_{\substack{s=1\\(s,s)\in \mathcal D}}^M {P^{(ss)}_{\mathrm{mic}}}^{\mathcal D}(\bG^{(ss)}),\\
&P^{\mathcal L}_{\mathrm{mic}}(\bG)
=\prod_{\substack{1\le s<t\le M\\(s,t)\in \mathcal L}}{P^{(st)}_{\mathrm{mic}}}^{\mathcal L}(\bG^{(st)})
\prod_{\substack{s=1\\(s,s)\in \mathcal L}}^M {P^{(ss)}_{\mathrm{mic}}}^{\mathcal L}(\bG^{(ss)}).
\end{aligned}
\ee
Here,
\begin{itemize}
\item
${P^{(st)}_{\mathrm{mic}}}^{\mathcal D}(\bG^{(st)})$ (${P^{(ss)}_{\mathrm{mic}}}^{\mathcal D}(\bG^{(ss)})$) 
is the microcanonical probability of the bi-partite (uni-partite) graph $\bG^{(st)}$ ($\bG^{(ss)}$) with 
constraints ${\vec{k}}^{\,*}_{s\to t}$ on the top layer and ${\vec{k}}^{\,*}_{t\to s}$ on the bottom layer 
(with constraint ${\vec{k}}^{\,*}_{s\to s}$).
\item
${P^{(st)}_{\mathrm{mic}}}^{\mathcal L}(\bG^{(st)})$ (${P^{(ss)}_{\mathrm{mic}}}^{\mathcal L}(\bG^{(ss)})$) 
is the microcanonical probability of the bi-partite (uni-partite) graph $\bG^{(st)}$ ($\bG^{(ss)}$) with 
constraint $L_{s,t}^*$ ($L^*_{s,s}$). 
\end{itemize}

The relative entropy becomes
\begin{eqnarray}
S_n(\Pmic \mid \Pcan)=S_n(P^{\mathcal D}_{\mathrm{mic}} \mid P^{\mathcal D}_{\mathrm{can}})
+S_n(P^{\mathcal L}_{\mathrm{mic}} \mid P^{\mathcal L}_{\mathrm{can}}).
\end{eqnarray}
It follows that
\be
\begin{aligned}
&\lim_{n_1,\dots,n_{M}\to\infty} \frac{S_n(\Pmic \mid \Pcan)}{n}\\
&\qquad =\lim_{n_1,\dots,n_{M}\to\infty}
\frac{S_n(P^{\mathcal D}_{\mathrm{mic}} \mid P^{\mathcal D}_{\mathrm{can}})}{n}
+\lim_{n_1,\dots,n_{M}\to\infty}
\frac{S_n(P^{\mathcal L}_{\mathrm{mic}} \mid P^{\mathcal L}_{\mathrm{can}})}{n}.
\end{aligned}
\ee
Using Theorem \ref{resultcomm} we get
\be
\begin{aligned}
&\lim_{n_1,\dots,n_{M}\to\infty}
\frac{S_n(P^{\mathcal D}_{\mathrm{mic}} \mid P^{\mathcal D}_{\mathrm{can}})}{n}
&= \sum_{(s,t)\in{\cal{D}}} A_s\,\| f_{s\to t}\|_{\ell^1(g)}.
\end{aligned}
\ee
Moreover,
\be
\begin{aligned}
&\lim_{n_1,\dots,n_{M}\to\infty}
\frac{S_n(P^{\mathcal L}_{\mathrm{mic}} \mid P^{\mathcal L}_{\mathrm{can}})}{n}\\
&= \lim_{n_1,\dots,n_{M}\to\infty}\sum_{\substack{1\le s<t\le M\\(s,t)\in \mathcal L}}
\frac{S_n({P^{(st)}_{\mathrm{mic}}}^{\mathcal L} \mid {P^{(st)}_{\mathrm{can}}}^{\mathcal L})}{n}
+\lim_{n_1,\dots,n_{M}\to\infty}\sum_{\substack{s=1\\(s,s)\in \mathcal L}}^M
\frac{S_n({P^{(ss)}_{\mathrm{mic}}}^{\mathcal L} \mid {P^{(ss)}_{\mathrm{can}}}^{\mathcal L})}{n}.
\end{aligned}
\ee
Using Theorems \ref{biplinks} and \ref{resultbip}, we get
\be
\lim_{n_1,\dots,n_{M}\to\infty}
\frac{S_n({P^{(st)}_{\mathrm{mic}}}^{\mathcal L} \mid {P^{(st)}_{\mathrm{can}}}^{\mathcal L})}{n}
=\lim_{n_1,\dots,n_{M}\to\infty}
\frac{S_n({P^{(ss)}_{\mathrm{mic}}}^{\mathcal L} \mid {P^{(ss)}_{\mathrm{can}}}^{\mathcal L})}{n}=0,
\ee
which proves the claim.
\end{proof}

\subsection{Proof of Theorem \ref{resultcommgen}}

\begin{proof}
The proof is based on the previous theorems. For each pair of layers $s,t\in {{\cal{M}}_1}$ we define 
Lagrange multipliers $\vt_{s\to t}=(\theta_{1}^t,\dots,\theta_{{n_s}}^t)$ and $\vt_{t\to s}=(\theta_{1}^s,
\dots,\theta_{{n_t}}^s)$. For each pair of layers $s,\in {{\cal{M}}_1},\  t\in {{\cal{M}}_2}$ we define 
$\vt_{s\to t}=(\theta_{1}^t,\dots,\theta_{{n_s}}^t)$. The Hamiltonian is 
\be
\begin{aligned}
&H\big(\bG \mid \vt_{s\to t};\ s \in {{\cal{M}}_1},\  t \in {{\cal{M}}_1} \cup {{\cal{M}}_2},\ \gamma_{s,t}(\bGamma)=1\big)\\
&= \sum_{\substack{s,t\in {{\cal{M}}_1}\\ \gamma_{s,t}(\bGamma)=1}} \vt_{s\to t}\vec{s}_{s\to t}(\bG)
+\sum_{\substack{s\in {{\cal{M}}_1}\\ \gamma_{s,s}(\bGamma)=1}}\vt_{s\to s}\vec{s}_{s\to s}(\bG)
+\sum_{\substack{s\in {{\cal{M}}_1},\  
t\in {{\cal{M}}_2}\\ \gamma_{s,t}(\bGamma)=1}}\vt_{s\to t}\vec{s}_{s\to t}(\bG)\\
&= H_{{{\cal{M}}_1}\to {{\cal{M}}_1}}+H_{{{\cal{M}}_1}\to {{\cal{M}}_2}},
\end{aligned}
\ee
with 
\be
\begin{aligned}
H_{{{\cal{M}}_1}\to {{\cal{M}}_1}} 
&= \sum_{\substack{s,t\in {{\cal{M}}_1}\\ \gamma_{s,t}(\bGamma)=1}}\vt_{s\to t}\vec{s}_{s\to t}(\bG)
+ \sum_{\substack{s\in {{\cal{M}}_1}\\ \gamma_{s,s}(\bGamma)=1}}\vt_{s\to s}\vec{s}_{s\to s}(\bG),\\
H_{{{\cal{M}}_1}\to {{\cal{M}}_2}}
&= \sum_{\substack{s\in {{\cal{M}}_1},\  
t\in {{\cal{M}}_2}\\ \gamma_{s,t}(\bGamma)=1}}\vt_{s\to t}\vec{s}_{s\to t}(\bG).
\end{aligned}
\ee
Consequently, the canonical ensemble is
\be
\Pcan(\bG)=P^{{{\cal{M}}_1}\to {{\cal{M}}_1}}_{\mathrm{can}}(\bG)
P^{{{\cal{M}}_1}\to {{\cal{M}}_2}}_{\mathrm{can}}(\bG)
\ee
with
\be
\begin{aligned}
P^{{{\cal{M}}_1}\to {{\cal{M}}_1}}_{\mathrm{can}}(\bG)
&= \prod_{\substack{s,t\in {{\cal{M}}_1}\\ 
\gamma_{s,t}(\bGamma)=1}}{P^{(st)}_{\mathrm{can}}}^{top,bot}(\bG^{(st)})
\prod_{\substack{s\in {{\cal{M}}_1}\\ 
\gamma_{s,s}(\bGamma)=1}} {P^{(ss)}_{\mathrm{can}}}^{}(\bG^{(ss)}),\\
P^{{{\cal{M}}_1}\to {{\cal{M}}_2}}_{\mathrm{can}}(\bG)
&= \prod_{\substack{s\in {{\cal{M}}_1},\  
t\in {{\cal{M}}_2}\\ \gamma_{s,t}(\bGamma)=1}}{P^{(st)}_{\mathrm{can}}}^{top}(\bG^{(st)}).
\end{aligned}
\ee
Here,
\begin{itemize}
\item
$\bG^{(st)}$ ($\bG^{(ss)}$) is the bi-partite (uni-partite) graph between layers $s$ and $t$ (itself) 
obtained from the multi-partite graph $\bG$. The $n_s\times n_t$ ($n_s\times n_s$) matrix 
representing this bi-partite (uni-partite) graph has, for each $i\in \Lambda_s$ and $j\in \Lambda_t$ 
(for each $i,j\in s$), elements $g_{i,j}(\bG^{(st)})=g_{i,j}(\bG)$ ($g_{i,j}(\bG^{(ss)})=g_{i,j}(\bG)$).
\item
${P^{(st)}_{\mathrm{can}}}^{top,bot}(\bG^{(st)})$ is the canonical probability of the bi-partite 
graph $\bG^{(st)}$ with constraints ${\vec{k}}^{\,*}_{s\to t}$ on the top layer and ${\vec{k}}^{\,*}_{t\to s}$ 
on the bottom layer.
\item 
${P^{(ss)}_{\mathrm{can}}}^{}(\bG^{(ss)})$ is the canonical probability of the uni-partite graph 
$\bG^{(ss)}$ with constraint ${\vec{k}}^{\,*}_{s\to s}$.
\item
${P^{(st)}_{\mathrm{can}}}^{top}(\bG^{(st)})$ is the canonical probability of the bi-partite graph 
$\bG^{(st)}$ with constraint ${\vec{k}}^{\,*}_{s\to t}$ on the top layer.
\end{itemize}

We can split the microcanonical probability as products of microcanonical probabilities for 
simpler cases. The number of graphs with such a type of constraints is 
\begin{eqnarray}
\Omega_{{\vec{k}}^{\,*}_{s\to t}; \ s \in {{\cal{M}}_1},\  
t \in {{\cal{M}}_1}\cup {{\cal{M}}_2},\ \gamma_{s,t}(\bGamma)=1}
=\Omega_{{\vec{k}}^{\,*}_{s\to t};\, s,t \in {{\cal{M}}_1},\ \gamma_{s,t}(\bGamma)=1}
\Omega_{{\vec{k}}^{\,*}_{s\to t};\,s \in {{\cal{M}}_1},\ t \in {{\cal{M}}_2},\ \gamma_{s,t}(\bGamma)=1}.
\end{eqnarray}
This means that the microcanonical probability can be factorised as
\begin{eqnarray}
\Pmic(\bG)=P^{{{\cal{M}}_1}\to {{\cal{M}}_1}}_{\mathrm{mic}}(\bG)P^{{{\cal{M}}_1}
\to {{\cal{M}}_2}}_{\mathrm{mic}}(\bG)
\end{eqnarray}
with 
\be
\begin{aligned}
P^{{{\cal{M}}_1}\to {{\cal{M}}_1}}_{\mathrm{mic}}(\bG)
&= \prod_{\substack{s,t\in {{\cal{M}}_1}\\ \gamma_{s,t}(\bGamma)=1}}
{P^{(st)}_{\mathrm{mic}}}^{top,bot}(\bG^{(st)})
\prod_{\substack{s\in {{\cal{M}}_1}\\ \gamma_{s,s}(\bGamma)=1}}
{P^{(ss)}_{\mathrm{mic}}}^{}(\bG^{(ss)}),\\
P^{{{\cal{M}}_1}\to {{\cal{M}}_2}}_{\mathrm{mic}}(\bG)
&= \prod_{\substack{s\in {{\cal{M}}_1},\  
t\in {{\cal{M}}_2}\\ \gamma_{s,t}(\bGamma)=1}}{P^{(st)}_{\mathrm{mic}}}^{top}(\bG^{(st)}).
\end{aligned}
\ee
Here,
\begin{itemize}
\item
${P^{(st)}_{\mathrm{mic}}}^{top,bot}(\bG^{(st)})$ is the microcanonical probability of the 
bi-partite graph $\bG^{(st)}$ with constraints ${\vec{k}}^{\,*}_{s\to t}$ on the top layer and 
${\vec{k}}^{\,*}_{t\to s}$ on the bottom layer. 
\item
${P^{(ss)}_{\mathrm{mic}}}^{}(\bG^{(ss)})$ is the microcanonical probability of the 
uni-partite graph $\bG^{(ss)}$ with constraint ${\vec{k}}^{\,*}_{s\to s}$.
\item
${P^{(st)}_{\mathrm{mic}}}^{top}(\bG^{(st)})$ is the microcanonical probability of the 
bi-partite graph $\bG^{(st)}$ with constraint ${\vec{k}}^{\,*}_{s\to t}$ on the top layer.
\end{itemize}

The relative entropy becomes
\begin{eqnarray}
S_n(\Pmic \mid \Pcan) = S_n(P^{{{\cal{M}}_1}\to {{\cal{M}}_1}}_{\mathrm{mic}} \mid 
P^{{{\cal{M}}_1}\to {{\cal{M}}_1}}_{\mathrm{can}})
+ S_n(P^{{{\cal{M}}_1}\to {{\cal{M}}_2}}_{\mathrm{mic}} \mid 
P^{{{\cal{M}}_1}\to {{\cal{M}}_2}}_{\mathrm{can}}).
\end{eqnarray}
It follows that
\be
\begin{aligned}
&\lim_{n_1,\dots,n_{M}\to\infty}\frac{S_n(\Pcan \mid \Pcan)}{n}\\
&=\lim_{n_1,\dots,n_{M}\to\infty}
\frac{S_n(P^{{{\cal{M}}_1}\to {{\cal{M}}_1}}_{\mathrm{mic}} 
\mid P^{{{\cal{M}}_1}\to {{\cal{M}}_1}}_{\mathrm{can}})}{n}
+\lim_{n_1,\dots,n_{M}\to\infty}
\frac{S_n(P^{{{\cal{M}}_1}\to {{\cal{M}}_2}}_{\mathrm{mic}} 
\mid P^{{{\cal{M}}_1}\to {{\cal{M}}_2}}_{\mathrm{can}})}{n}.
\end{aligned}
\ee
Using again Theorem \ref{resultcomm} we get
\be
\begin{aligned}
&\lim_{n_1,\dots,n_{M}\to\infty}
\frac{S_n(P^{{{\cal{M}}_1}\to {{\cal{M}}_1}}_{\mathrm{mic}} 
\mid P^{{{\cal{M}}_1}\to {{\cal{M}}_1}}_{\mathrm{can}})}{n}\\
&= \sum_{\substack{s,t\in {{\cal{M}}_1}\\ 
\gamma_{s,t}(\bGamma)=1}}\left\lbrace A_s\,\|f_{s\to t}\|_{\ell^1(g)}
+A_t\,\|f_{t\to s}\|_{\ell^1(g)}\right\rbrace 
+ \sum_{\substack{s\in {{\cal{M}}_1}\\ \gamma_{s,s}(\bGamma)=1}} 
A_s\,\|f_{s\to s}\|_{\ell^1(g)}\\
&=\sum_{\substack{s, t\in {{\cal{M}}_1}\\ \gamma_{s,t}(\bGamma)=1}}
 A_s\,\|f_{s\to t}\|_{\ell^1(g)}.
\end{aligned}
\ee
From Theorem~\ref{configurationmodel} we get
\be
\begin{aligned}
&\lim_{n_1,\dots,n_{{M}}\to\infty}
\frac{S_n(P^{{{\cal{M}}_1}\to {{\cal{M}}_2}}_{\mathrm{mic}} 
\mid P^{{{\cal{M}}_1}\to {{\cal{M}}_2}}_{\mathrm{can}})}{n}\\
&= \lim_{n_1,\dots,n_{{M}}\to\infty}
\sum_{\substack{s\in {{\cal{M}}_1},\  
t\in {{\cal{M}}_2}\\ \gamma_{s,t}(\bGamma)=1}}
\frac{S_n({P^{(st)}_{\mathrm{mic}}}^{{{{\cal{M}}_1}\to {{\cal{M}}_2}}} \mid 
{P^{(st)}_{\mathrm{can}}}^{{{{\cal{M}}_1}\to {{\cal{M}}_2}}})}{n}
=\sum_{\substack{s\in {{\cal{M}}_1},\  
t\in {{\cal{M}}_2}\\ \gamma_{s,t}(\bGamma)=1}} 
A_s\,\|f_{s\to l}\|_{\ell^1(g)},
\end{aligned}
\ee
which concludes the proof.
\end{proof}



\begin{thebibliography}{99}

\bibitem{ginestra}
K.\ Anand and G.\ Bianconi,
\emph{Phys.\ Rev.\ E} \textbf{80}, 045102 (2009).  

\bibitem{gin_cavity}
K.\ Anand and G.\ Bianconi,
\emph{Phys.\ Rev.\ E} \textbf{82}, 011116 (2010).

\bibitem{arratialiggett2005}
R.\ Arratia and T.M.\ Liggett, 
\emph{Ann.\ Appl.\ Probab.} \textbf{15}, 652 (2005).

\bibitem{Barre2007}
J.\ Barr\'e and B.\ Goncalves, 
\emph{Physica A} \textbf{386}, 212 (2007).

\bibitem{Barre2001}
J.\ Barr\'e, D.\ Mukamel, and S.\ Ruffo, 
\emph{Phys.\ Rev.\ Lett.} \textbf{87}, 030601 (2001).

\bibitem{bender1977}
E.A.\ Bender, 
\emph{Discrete Math.} \textbf{10}, 217 (1974).

\bibitem{bender1978}
E.A.\ Bender and E.R.\ Canfield, 
\emph{J.\ Comb.\ Theory (A)} \textbf{24}, 296 (1978).

\bibitem{bianconiMultiplex}
G.\ Bianconi, 
\emph{Physical Review E} \textbf{87}, 062806 (2013).

\bibitem{gin_hierarchically}
G.\ Bianconi, A.C.C.\ Coolen, C.J.\ Perez Vicente,
\emph{Phys.\ Rev.\ E} \textbf{78}, 016114 (2008).

\bibitem{Blume1971}
M.\ Blume, V.J.\ Emery, and R.B.\ Griffiths, 
\emph{Phys.\ Rev.\ A} \textbf{4}, 1071 (1971).

\bibitem{multiplex}
S.\ Boccaletti, G.\ Bianconi, R.\ Criado, C.I.\ Del Genio, J.\ G\'omez-Garde\~{n}es, 
M.\ Romance, I.\ Sendi\~{n}a-Nadalj, Z.\ Wang, and M.\ Zanin, 
\emph{Physics Reports} \textbf{544}, 1 (2014).

\bibitem{cutoff}
M.\ Bogu\~{n}\'a, R.\ Pastor-Satorras, and A.\ Vespignani, 
\emph{The European Physical Journal B} \textbf{38}, 205-209 (2004).

\bibitem{bollobas1980}
B.\ Bollob\'as, 
\emph{European J.\ Combin.} \textbf{1}, 311 (1980).

\bibitem{guidosbook}
G.\ Caldarelli, \emph{Scale-Free Networks: Complex Webs in Nature and Technology} 
(OUP, Oxford, 2007).

\bibitem{Ruffo}
A.\ Campa, T.\ Dauxois, D.\ Fanelli, and S.\ Ruffo, 
\emph{Physics of Long-Range Interacting Systems} 
(OUP, Oxford, 2014).

\bibitem{Ruffo2}
A.\ Campa, T.\ Dauxois, and S.\ Ruffo, 
\emph{Physics Reports} \textbf{480}, 57 (2009).

\bibitem{WF}
P.J.\ Carrington, J.\ Scott, and S.\ Wasserman, 
\emph{Models and Methods in Social Network Analysis}, 
Vol.\ 28 (Cambridge University Press, 2005).

\bibitem{chatterjee}
S.\ Chatterjee, and P.\ Diaconis,
\emph{The Annals of Statistics} \textbf{41}, 2428 (2013).

\bibitem{Chavanis2003}
P.-H.\ Chavanis, 
\emph{Astron.\ \& Astrophys.} \textbf{401}, 15 (2003).

\bibitem{chunglu}
F.\ Chung, and L.\ Lu, 
\emph{Proceedings of the National Academy of Sciences} \textbf{99}, 15879 (2002).

\bibitem{DAgostino2000}
M.\ D'Agostino \emph{et al.}, 
\emph{Phys.\ Lett.\ B} \textbf{473}, 219 (2000).

\bibitem{netonets}
G.\ D'Agostino, and A.\ Scala, 
\emph{Networks of Networks: The Last Frontier of Complexity}, 
Vol.\ 340 (Springer, Berlin, 2014).

\bibitem{Ellis2000}
R.S.\ Ellis, K.\ Haven, and B.\ Turkington, 
\emph{J.\ Stat.\ Phys.} \textbf{101}, 999 (2000).

\bibitem{Ellis2002}
R.S.\ Ellis, K.\ Haven, and B.\ Turkington, 
\emph{Nonlinearity} \textbf{15}, 239 (2002).

\bibitem{Ellis2004}
R.S.\ Ellis, H.\ Touchette, and B.\ Turkington, 
\emph{Physica A} \textbf{335}, 518 (2004).

\bibitem{community_santo}
S.\ Fortunato, 
\emph{Physics Reports} \textbf{486}, 75 (2010).

\bibitem{fronczakBM}
P.\ Fronczak, A.\ Fronczak, and M.\ Bujok, 
\emph{Physical Review E} \textbf{88}, 032810 (2013).

\bibitem{myX}
D.\ Garlaschelli, 
\emph{New Journal of Physics} \textbf{11}, 073005 (2009).

\bibitem{myY}
D.\ Garlaschelli, S.E.\ Ahnert, T.\ Fink, and G.\ Caldarelli, 
\emph{Entropy} \textbf{15}, 3148 (2013).

\bibitem{GdHR15}
D.\ Garlaschelli, F.\ den Hollander, and A.\ Roccaverde,
\emph{Nieuw Archief voor Wiskunde} \textbf{5/16}, 207 (2015).

\bibitem{mylikelihood}
D.\ Garlaschelli and M.I.\ Loffredo, 
\emph{Physical Review E} \textbf{78}, 015101 (2008).

\bibitem{mymultiplexity}
V.\ Gemmetto, and D.\ Garlaschelli, 
\emph{Scientific Reports} \textbf{5}, 9120 (2015).

\bibitem{Gibbs1902}
J.W.\ Gibbs, 
\emph{Elementary Principles of Statistical Mechanics}, 
Yale University Press, New Haven, Connecticut (1902).

\bibitem{GMW06}
C.\ Greenhill, B.D.\ McKay, and X.\ Wang,
\emph{J.\ Comb.\ Theory, Series A} \textbf{113}, 291 (2006).

\bibitem{Thirring1971}
P.\ Hertel and W.\ Thirring, 
\emph{Annals of Physics} \textbf{63}, 520 (1971).

\bibitem{vdH16}
R.\ van der Hofstad, 
\emph{Random Graphs and Complex Networks, Volume I}, 
to appear with Cambridge University Press. 

\bibitem{simpleBM}
P.W.\ Holland, K.B.\ Laskey, and S.\ Leinhardt, 
\emph{Social Networks} \textbf{5}, 109 (1983).

\bibitem{timevarying}
P.\ Holme, and J.\ Saram\"{a}ki, 
\emph{Physics Reports} \textbf{519}, 97 (2012).

\bibitem{degreecorrectedBM}
B.\ Karrer, and M.E.J.\ Newman, 
\emph{Physical Review E} \textbf{83}, 016107 (2011).

\bibitem{LyndenBell1999}
D.\ Lynden-Bell, 
\emph{Physica A} \textbf{263}, 293 (1999).

\bibitem{LyndenBell1968}
D.\ Lynden-Bell, R.\ Wood, 
\emph{Monthly Notices of the Royal Astronomical Society} \textbf{138}, 495 (1968).

\bibitem{myenhanced}
R.\ Mastrandrea, T.\ Squartini, G.\ Fagiolo, and D.\ Garlaschelli, 
\emph{New Journal of Physics} \textbf{16}, 043022 (2014).

\bibitem{mckay1991}
B.D.\ McKay and N.C.\ Wormald, 
\emph{Combinatorica} \textbf{11}, 369 (1991).

\bibitem{molloyreed}
M.\ Molloy, and B.\ Reed, 
\emph{Random Structures \& Algorithms} \textbf{6}, 161 (1995).

\bibitem{newmanstrogatz}
M.E.\ Newman, S.H.\ Strogatz, and D.J.\ Watts, 
\emph{Physical review E} \textbf{64}, 026118 (2001).

\bibitem{gin_growmultiplex}
V.\ Nicosia, G.\ Bianconi, V.\ Latora, and M.\ Barthelemy,
\emph{Phys.\ Rev.\ Lett.} \textbf{111}, 058701 (2013).

\bibitem{ParkNewman}
J.\ Park, and M.E.\ Newman, 
\emph{Physical Review E} \textbf{70}, 066117 (2004).

\bibitem{entropyBM}
T.P.\ Peixoto, 
\emph{Physical Review E} \textbf{85}, 056122 (2012).

\bibitem{community_porter}
M.A.\ Porter, J.P.\ Onnela, and P.J.\ Mucha, 
\emph{Notices of the AMS} \textbf{56}, 1082-1097.

\bibitem{RS}
C.\ Radin, and L.\ Sadun, 
\emph{J.\ Phys.\ A:\ Math.\ Theor.} \textbf{46}, 305002  (2013).

\bibitem{RS3}
C.\ Radin, and L.\ Sadun, 
\emph{Journal of Statistical Physics} \textbf{158}, 853 (2015).

\bibitem{RS2}
C.\ Radin and M.\ Yi,
\emph{The Annals of Applied Probability} \textbf{23}, 2458-2471 (2013).

\bibitem{mymethod}
T.\ Squartini, and D.\ Garlaschelli, 
\emph{New Journal of Physics} \emph{13}, 083001 (2011).

\bibitem{myunbiased}
T.\ Squartini, R.\ Mastrandrea, and D.\ Garlaschelli, 
\emph{New J.\ Phys.} \textbf{17}, 023052  (2015).

\bibitem{SdMdHG15}
T.\ Squartini, J.\ de Mol, F.\ den Hollander, and D.\ Garlaschelli, 
\emph{Phys.\ Rev.\ Lett.} \textbf{115}, 268701 (2015).

\bibitem{Thirring1970}
W.\ Thirring, 
\emph{Zeitschrift f\"ur Physik} \textbf{235}, 339 (1970).

\bibitem{touchette2014}
H.\ Touchette, 
\emph{J. Stat. Phys.} \textbf{159}, 987 (2015).

\bibitem{Touchette2004}
H.\ Touchette, R.S.\ Ellis, and B.\ Turkington, 
\emph{Physica A} \textbf{340}, 138 (2004).

\bibitem{ziff}
R.M.\ Ziff, G.E.\ Uhlenbeck, and M.\ Kac,
\emph{Physics Reports} \textbf{32}, 169 (1977).

\end{thebibliography}
\end{document}